\documentclass[12pt,reqno]{amsproc}
\usepackage[all]{xy}
\usepackage{amsmath,graphicx,amssymb,amsthm, mathrsfs, amsfonts,dsfont}
\usepackage{rotating}
\usepackage{CJKutf8}
\allowdisplaybreaks
\def\A{\mathcal A}
\def\E{\mathcal E}

\def\U{\mathcal U}

\def\H{\mathcal H}

\def\M{\mathcal M}
\def\N{\mathcal N}

\def\Z{\mathcal Z}

\def\amslatex{$\mathcal{A}\kern-.1667em\lower.5ex\hbox{$\mathcal{M}$}\kern-.125em\mathcal{S}$-\LaTeX}

\newtheorem{set}{set}[section]
\newtheorem{Corollary}[set]{Corollary}

\newtheorem{Lemma}[set]{Lemma}

\newtheorem{Remark}[set]{Remark}
\newtheorem{Theorem}[set]{Theorem}
\newcommand{\define}{\mathrel{\hbox{$\equiv$\hskip -.90em \lower .47ex \hbox{$\leftharpoondown$}}}}
\newcommand{\enifed}{\mathrel{\hbox{$\equiv$\hskip -.90em \lower .47ex \hbox{$\rightharpoondown$}}}}

\numberwithin{equation}{section}
\usepackage{blkarray}
\makeatletter
\newcommand{\LeftEqNo}{\let\veqno\@@leqno}

\makeatother

\begin{document}
\begin{CJK}{UTF8}{<font>}
\title{A local quantization principle for inclusions of tracial von Neumann algebras}

\author{Xinyan Cao}
\address{School of Mathematical Sciences, Dalian University of Technology, Dalian, Liaoning, 116024, China}
\email{caoxinyan123@qq.com}

\author{Junsheng Fang}
\address{School of Mathematical Sciences, Hebei Normal University, Shijiazhuang, Hebei, 050024, China}
\email{jfang@hebtu.edu.cn}
\thanks{Junsheng Fang was partly supported by NSFC(Grant No.12071109), a Start-up funding of Hebei Normal University, and the Hebei Natural Science Foundation (No.
A2023205045).}

\author{Chunlan Jiang}
\address{School of Mathematical Sciences, Hebei Normal University, Shijiazhuang, Hebei, 050024, China}
\email{cljiang@hebtu.edu.cn}
\thanks{ Chunlan Jiang was partly supported by National Natural Science Foundation
of China (NO. 11831006), and the Hebei Natural Science Foundation (No.
A2023205045).}

\author{Zhaolin Yao}
\address{School of Mathematical Sciences, Hebei Normal University, Shijiazhuang, Hebei, 050024, China}
\email{zyao@hebtu.edu.cn}
\thanks{Zhaolin Yao was partly supported by NSFC (Grant No.12371128), Science Foundation
of Hebei Normal University (Grant No. L2023B05) and the Hebei Natural Science Foundation (No.
A2023205045).}
\date{}
\maketitle
\begin{abstract}
We study the local quantization principle (after Sorin Popa~\cite{popa 94} and \cite{popa 95}) of inclusions of tracial von Neumann algebras. Let $(\mathcal{M},\tau)$ be a type ${\rm II}_1$ von Neumann algebra and let $\mathcal{N}\subseteq \mathcal{M}$ be a type ${\rm II}_1$ von Neumann subalgebra. Let $x_1,\ldots, x_m \in \mathcal{M}$ and $ \epsilon> 0$. Then there exists a partition of 1 with projections $p_{1}, \ldots, p_{n}$ in $\mathcal{N}$ such that
\[\left\|\sum_{i=1}^n p_{i}\left(x_j-E_{\mathcal{N}'\cap \mathcal{M}}(x_j)\right)p_{i}\right\|_{2}<\epsilon,\quad 1\leq j\leq m.\] In particular, if $\N\subseteq \M$ is an inclusion of type $\rm II_{1}$ factors with $[\M:\N]=2$, then for any $x_{1},\ldots, x_{m}\in \M$, there exists a partition of 1 with projections $p_{1}, \ldots, p_{n}$ in $\mathcal{N}$ such that \[\sum_{i=1}^n p_ix_jp_i=\tau(x_j)1, \quad 1\leq j\leq m.\] Equivalently, there exists a unitary operator $u\in \N$ such that \[\frac{1}{n}\sum_{i=1}^nu^{*i}x_j u^i=\tau(x_j)1, \quad 1\leq j\leq m.\]

\end{abstract}

\vskip1.0cm

\section{Introduction}

In von Neumann algebras, the Dixmier property (see \cite{kR} Theorem 8.3.6) is one of the fundamental approximation properties. This important result states that for any element $x$ in a von Neumann algebra $\mathcal{M}$ and for any $\epsilon>0$, there are unitary operators $u_1, \ldots,u_n$ in $\mathcal{M}$ and positive numbers $\lambda_1,\ldots,\lambda_n$ of sum 1 such that
\[\left\|\sum_{i=1}^n \lambda_i u_i^*x u_i-E_{\Z(\mathcal{M})}(x)\right\| <\epsilon.
\]
A natural question arises: if $\mathcal{N}\subset \mathcal{M}$ is a von Neumann subalgebra, is above property still true
under the additional requirement that the unitaries $u_i$ belong to $\mathcal{N}$ and $\Z(\mathcal{M})$ is replaced with $\mathcal{N}^\prime \cap \mathcal{M}$? If the answer is yes, then we say that the inclusion $\mathcal{N}\subset \mathcal{M}$ has relative Dixmier property. The relative Dixmier property has been studied by several authors (see \cite{RLA,HVG,Jaeseong,pop,Popa 2000}), in connection with extensions of pure states, a problem originating in \cite{kR59}.

In this paper, we focus on a weaker version of the Dixmier property which is more suitable for von Neumann algebras. Let $\mathcal{N}\subset \mathcal{M}$ be an inclusion of tracial von Neumann algebras. Following \cite{Popa 2000} Definition 1.1, we say that the inclusion $\mathcal{N} \subset \mathcal{M}$ has the weak relative Dixmier property if for every $x\in \mathcal{M}$, the weak* closed convex hull of $\{uxu^* | u\in \U(\mathcal{N})\}$ intersects $\mathcal{N}^\prime \cap \mathcal{M}$. With the above notation, if there are mutually orthogonal projections $p_1,\ldots,p_n$ in $\mathcal{N}$ of sum 1 such that
\[\left\|\sum_{i=1}^n p_ixp_i-E_{\mathcal{N}^\prime \cap \mathcal{M}}(x)\right\|_2 <\epsilon,
\]
then the weak relative Dixmier property follows easily from above inequality. The type of result that we will discuss proves to be very useful in exploiting the noncommutative ergodic phenomena specific to the theory of type $\rm II_1$ factors (see \cite{popa 81, popa 83, popa 86, popa 90, popa 88}). Sorin Popa proved the following general result (see Lemma A.1.1 of~\cite{popa 94}), which plays a key role in this paper.

\begin{Theorem}
Let $\mathcal{N}\subseteq \mathcal{M}$ be tracial von Neumann algebras with a normal finite faithful trace $\tau$. Let $X\subset\mathcal{M}$ be a finite set of elements such that $E_{\mathcal{N}\vee(\mathcal{N}'\cap \mathcal{M})}(x)=0$ for all $x\in X$. Given any $ \epsilon> 0$, there exists a partition of the identity with projections $\{p_i\}_{i=1}^n$ in $\mathcal{N}$ such that $\|\sum_{i=1}^n p_ixp_i\|_{2}<\epsilon$ for all $x\in X $.
\end{Theorem}

As a quick corollary, let $\N$ is a maximal abelian von Neumann subalgebra of a tracial von Neumann algebra $(\M,\tau)$ and $x_1,\ldots,x_m\in \M$ with $E_{\mathcal{N}}(x_i)=0$ for $i=1,\ldots,m$. Given any $\epsilon>0$, there exists a partition of the identity with projections $\{p_i\}_{i=1}^n$ in $\N$ such that $\|\sum_{i=1}^np_ix_j p_i\|_{2}<\epsilon$ for all $1\leq j\leq m$. In this paper, we prove the following main theorem.

\begin{Theorem}\label{T:local quantization principle}
Let $(\mathcal{M},\tau)$ be a type ${\rm II}_1$ von Neumann algebra and let $\mathcal{N}\subseteq \mathcal{M}$ be a type ${\rm II}_1$ von Neumann subalgebra. Let $x_1,\ldots,x_m\in \mathcal{M}$ be a finite set of elements and $ \epsilon> 0$. Then there exists a partition of 1 with projections $p_{1}, \ldots, p_{n}$ in $\mathcal{N}$ such that
\[\left\|\sum_{i=1}^n p_{i}\left(x_j-E_{\mathcal{N}'\cap \mathcal{M}}(x_j)\right)p_{i}\right\|_{2}<\epsilon,\quad \forall\, 1\leq j\leq m.\]
\end{Theorem}

Furthermore, if we assume that $\N\subseteq \M$ are type $\rm II_{1}$ factors and $[\M:\N]=2$, we will prove the following result.

\begin{Theorem}\label{index 2 finite proj in N}
Let $(\mathcal{M},\tau)$ be a type ${\rm II}_1$ factor and $\N\subseteq \M$ be an inclusion of type $\rm II_{1}$ subfactor with $[\M:\N]=2$. For any $x_{1},\ldots, x_{m}\in \M$, there exists a partition of 1 with projections $p_{1},\ldots, p_{n}$ in $\mathcal{N}$ such that $\sum_{i=1}^np_ix_jp_i=\tau(x_j)1$ for $1\leq j\leq m$.
\end{Theorem}

The organization of this paper is as follows. In section 2, we prove Theorem 1.2 stated as above. In section 3, we prove Theorem 1.3. We refer to ~\cite{kR,SZ} for the basic knowledge of von Neumann algebras. We refer to~\cite{AP,SS} for the general theory of type ${\rm II}_1$ factors.

\section{Proof of Theorem~\ref{T:local quantization principle}}
To prove Theorem~\ref{T:local quantization principle}, we need the following lemmas.

\begin{Lemma}\label{11}
Let $\mathcal{M}$ be a von Neumann algebra and $x\in\mathcal{M}$. Then there is a nonzero projection $e\in \mathcal{M}$ such that $exe=0$ if and only if there is a nonzero element $y\in\mathcal{M}$ such that $yxy^{*}=0$.
\end{Lemma}
\begin{proof}
Suppose $y\neq 0$ satisfies $yxy^*=0$. By the polar decomposition theorem, there is a nonzero positive operator $h=(y^*y)^{1/2}\in\M$ and a partial isometry $v\in\M$ such that $y=vh$ and $v^*vh=h=hv^*v$. Then $yxy^*=0$ implies that
$
vhxhv^*=0
$
and therefore,
\[
hxh=v^*vhxhv^*v=0.
\]
Let $f=\chi_{h}(\{0\})$ and $e=I-f$. Then $e\neq 0$ and
\[
h=\begin{pmatrix}
h_{1}&0\\
0&0
\end{pmatrix}
\]
with respect to the decomposition $I=e+f$, where $h_{1}=ehe$ is an injective positive operator with dense range.
Write
\[
x=\begin{pmatrix}
x_{11}&x_{12}\\
x_{21}&x_{22}
\end{pmatrix}
\]
with respect to the decomposition $I=e+f$. Then
\[
0=hxh=\begin{pmatrix}
h_{1}&0\\
0&0
\end{pmatrix}
\begin{pmatrix}
x_{11}&x_{12}\\
x_{21}&x_{22}
\end{pmatrix}
\begin{pmatrix}
h_{1}&0\\
0&0
\end{pmatrix}=
\begin{pmatrix}
h_{1}x_{11}h_1&0\\
0&0
\end{pmatrix}
\]
implies that $h_1x_{11}h_1=0$. Since $h_1$ is injective with dense range, it follows that $x_{11}=0$. This implies that $exe=0$.

\end{proof}

We denote by $\Z(\M)$ the center of the von Neumann algebra $\M$, that is, $\Z(\M)=\M\cap\M'$.

\begin{Lemma}\label{10}
Let $\mathcal{M}$ be a von Neumann algebra and let $u\in\mathcal{M}$ be a unitary operator. If $u=u^{*}$ and $u\notin \mathcal{Z}(\M)$, then there is a nonzero projection $q\in \mathcal{M}$ such that $quq=0$.
\end{Lemma}
\begin{proof}
Let $f_{1}=\chi_{u}(\{1\})$ and $f_{2}=\chi_{u}(\{-1\})$. Then $f_1-f_2=u$ and $f_{1}+f_{2}=I$. By the Comparison theorem, there exists a central projection $p$ such that $pf_{1}\preceq pf_{2}$ and $(I-p)f_{1}\succeq(I-p)f_{2}$. Since $u\notin \mathcal{Z}(\M)$,  it follows that $f_1, f_2\notin \mathcal{Z}(\M)$ and so at least one of $pf_1$ or $(I-p)f_2$ must be nonzero. We may assume that $0<pf_{1}\preceq pf_{2}$, then there is a subprojection $G$ of $pf_{2}$ such that $pf_{1}\sim G$ in $\mathcal{M}$. Note that \[(pf_{1}+G)\mathcal{M}(pf_{1}+G)\cong M_{2}(\mathbb{C})\otimes(pf_{1}\mathcal{M}pf_{1}).\] Let $v=pf_{1}-G$. Then we can view $v$ in $(pf_1+G)\M(pf_1+G)$ as the element $\begin{pmatrix}
  1&0 \\
0&-1
\end{pmatrix}$ in $M_2(\mathbb{C})\otimes (pf_1\M pf_1)$. A simple computation shows that
\[
\begin{pmatrix}
  \frac{1}{2}&\frac{1}{2} \\
 \frac{1}{2}&\frac{1}{2}
\end{pmatrix}
\begin{pmatrix}
  1&0 \\
0&-1
\end{pmatrix}
\begin{pmatrix}
  \frac{1}{2}&\frac{1}{2} \\
 \frac{1}{2}&\frac{1}{2}
\end{pmatrix}=0.
\]
Therefore, there is a nonzero subprojection $q$ of $pf_{1}+G$ such that $qvq=0$. Thus \[quq=q(pf_{1}+G)u(pf_{1}+G)q=q(pf_{1}+G)(f_{1}-f_{2})(pf_{1}+G)q=qvq=0.\]

\end{proof}

\begin{Lemma}\label{9}
Let $\mathcal{M}$ be a finite von Neumann algebra. For any self-adjoint operator $x\in\mathcal{M}$, let $x=uh$, where $u$ is a unitary operator. If $u\not \in \mathcal{Z}(\M)$, then there is a nonzero projection $p\in \mathcal{M}$ such that $pxp=0$.
\end{Lemma}
\begin{proof}
If $\chi_{x}(\{0\})\neq0$, then let $p=\chi_{x}(\{0\})$. We have $p\neq 0$ and $pxp=0$. So we may assume that $\chi_{x}(\{0\})=0$, i.e., $x$ is injective with dense range. Since $u\notin \Z(\mathcal{M})$, by lemma \ref{10}, there is a nonzero projection $f\in \mathcal{M}$ such that $fuf=0$.
From this and $uh=hu$, we have $fh^{-\frac{1}{2}}uhh^{-\frac{1}{2}}f=0$ and $fh^{-\frac{1}{2}}\neq0$. Let \[z=fh^{-\frac{1}{2}}.\] Then $z$ is a closed densely defined operator affiliated with $\mathcal{M}$ and $z^{*}=h^{-\frac{1}{2}}f$. Let \[y=(I+zz^{*})^{-\frac{1}{2}}z.\] Then $y$ is a bounded operator in $\mathcal{M}$ and $y\neq0$. Note that $y^{*}=z^{*}(I+zz^{*})^{-\frac{1}{2}}$ and $yxy^{*}=0$. By lemma \ref{11}, there is a nonzero projection $p\in \mathcal{M}$ such that $pxp=0$.

\end{proof}

\begin{Lemma}\label{EXE=0}
Let $\M$ be a finite von Neuamnn algebra and let $x\in\M$ be a self-adjoint operator. If $E_{\mathcal{Z}(\M)}(x)=0$, then there is a nonzero projection $p$ in $\M$ such that $pxp=0$.
\end{Lemma}
\begin{proof}
If $\text{ker} x\neq \{0\}$, then let $p$ be the projection onto $\text{ker} x$. We have $p\neq 0$ and $pxp=0$. Otherwise, $\text{ker} x=0$ and the range projection of $x$ is $I$. Let $x=uh$ be the polar decomposition of $x$. Then $u$ is a unitary operator and $h\geq 0$ is injective with dense range. If $u\in \mathcal{Z}(\M)$, then \[0=E_{\mathcal{Z}(\M)}(uh)=uE_{\mathcal{Z}(\M)}(h).\] So $E_{\mathcal{Z}(\M)}(h)=0$ and $h=0$. This is a contradiction. Thus $u\notin \mathcal{Z}(\M)$. By lemma \ref{9}, there is a nonzero projection $p\in \M$ such that $pxp=0$.

\end{proof}

\begin{Corollary}\label{8}
Let $(\mathcal{M}, \tau)$ be a tracial von Neumann algebra and let $x\in\M$ be a self-adjoint operator. If $E_{\mathcal{Z}(\M)}(x)=0$, then there is a sequence of nonzero projections $\{p_n\}_{n=1}^\infty$ in $\M$ such that $\sum_{n=1}^\infty p_n=I$ and $p_nxp_n=0$ for all $n\geq 1$.
\end{Corollary}
\begin{proof}
Let $\E=\{\{e_\alpha\}\}$ be the family of mutually orthogonal nonzero projections such that $e_\alpha xe_\alpha=0$ for all $\alpha$. By lemma \ref{EXE=0}, $\E$ is not empty. Define $\{e_\alpha\}\prec \{f_\beta\}$ if $\{e_\alpha\}\subseteq \{f_\beta\}$. Clearly, $\prec$ is a partial order and if $\{\{e_\alpha\}\}$ is a totally ordered set, then the union is a maximal element of the totally ordered set. By Zorn's lemma, there is a maximal element $\{e_\alpha\}$ in $\E$. We claim that $\sum_\alpha e_\alpha=I$. Otherwise, let $f=I-\sum_{\alpha}e_\alpha$.
Then
\begin{align*}
 E_{\Z(\M)}(fxf)&=E_{\Z(\M)}(xf)=E_{\Z(\M)}\left(x-x\sum_\alpha e_\alpha\right)\\
               &=E_{\Z(\M)}(x)-\sum_\alpha E_{Z(\M)}(e_\alpha xe_\alpha)=0.
\end{align*}
We consider $y=fxf$ as an element in $f\M f$. If ker$y\neq \{0\}$. Let $e$ be the projection onto ker$y$. We have $0\neq e\leq f$ and $exe=efxfe=eye=0$. Assume that $\text{ker}y=\{0\}$ and $y=uh$ is the polar decomposition of $y$. Then $u$ is a unitary operator in $f\M f$ and $0\neq h\in f\M f$. We claim that $u$ is not in $\Z(f\M f)$. Otherwise, $u\in \Z(f\M f)$. Write $u=e^{ia}$ and $a=a^*\in \Z(f\M f)$. Note that \[f\M f\cap f\M' f=\Z(\M)f.\] Therefore, there is an operator $b=b^*\in \Z(\M)$ such that $a=bf$. Therefore, $u=e^{ib}f$. Since $0\neq h\in f\M f$, we have
\[
0=E_{\Z(\M)}(y)=E_{\Z(\M)}(uh)=E_{\Z(\M)}(e^{ib}fh)=E_{\Z(\M)}(e^{ib}h)=e^{ib}E_{\Z(\M)}(h).
\]
This implies that $E_{\Z(\M)}(h)=0$ and so $h=0$, which is a contradiction. By lemma~\ref{9}, there is a nonzero projection $e\leq f$ such that \[efxfe=0.\] Thus $exe=0$ and $\{e_\alpha\}\cup \{e\}\in \E$. This is a contradiction.

\end{proof}

The following lemma is motivated by Lemma 2.3 of \cite{WC} and its proof is similar to the proof of Lemma 2.3 of \cite{WC}. We provide a complete proof for the reader¡¯s convenience.

\begin{Lemma}\label{Y1-1Y^{*}}
Let $\M$ be a finite von Neumann algebra and $a, b\in \M$. Then there is a nonzero operator $y\in M_{2}(\M)$ such that
\[y
\begin{pmatrix}
1& a\\
b&-1
\end{pmatrix}
y^{*}=0.
\]
\end{Lemma}
\begin{proof}
By the polar decomposition theorem, there is a unitary operator $u\in \M$ and a positive operator $h\in \M$ such that $a^{*}-b=uh$. Then \[u^{*}(a^{*}-b)=u^{*}uh=h\] is a positive operator. Note that
\[u^{*}(a^{*}-b)=(u^{*}(a^{*}-b))^{*}=(a-b^{*})u.\]
Then \[(u^{*}b+au)^{*}=b^{*}u+u^{*}a^{*}=u^{*}b+au.\]  Therefore, $u^{*}b+au$ is a self-adjoint operator. Since $E_{\mathcal{Z}(\M)}(u^{*}b+au)\in\mathcal{Z}(\M)$ and $\mathcal{Z}(\M)$ is abelian, we can write \[\phi(\omega)=E_{\mathcal{Z}(\M)}(u^{*}b+au).\] Let \[t(\omega)=\frac{\phi(\omega)+\sqrt{\phi(\omega)^{2}+4}}{2}\in \mathcal{Z}(\M).\]
A simple computation shows that
\begin{eqnarray}\nonumber
~&&\begin{pmatrix}
1& t(w)u^{*}\\
1&t(w)u^{*}
\end{pmatrix}
\begin{pmatrix}
1& a\\
b&-1
\end{pmatrix}
\begin{pmatrix}
1& 1\\
t(w)u&t(w)u
\end{pmatrix}  \\\nonumber
   &=&
\begin{pmatrix}
1+t(w)(u^{*}b+au)-t(w)^{2}& 1+t(w)(u^{*}b+au)-t(w)^{2}\\
1+t(w)(u^{*}b+au)-t(w)^{2}&1+t(w)(u^{*}b+au)-t(w)^{2}
\end{pmatrix}.
\end{eqnarray}
Note that \[E_{\mathcal{Z}(\M)}(1+t(w)(u^{*}b+au)-t(w)^{2})=1+t(w)\phi(w)-t(w)^2=0.\] By lemma \ref{EXE=0}, there is a nonzero projection $e\in \M$ such that
\[e(1+t(w)(u^{*}b+au)-t(w)^{2})e=0.\]
Let
\[y=
\begin{pmatrix}
e& 0\\
0&0
\end{pmatrix}
\begin{pmatrix}
1& t(w)u^{*}\\
1&t(w)u^{*}
\end{pmatrix}=
\begin{pmatrix}
e& t(w)eu^{*}\\
0&0
\end{pmatrix}\neq0.
\]
Thus
\begin{eqnarray}\nonumber
~&& y
\begin{pmatrix}
1& a\\
b&-1
\end{pmatrix}
y^{*} \\\nonumber
&=&
\begin{pmatrix}
e& 0\\
0&0
\end{pmatrix}
\begin{pmatrix}
1& t(w)u^{*}\\
1&t(w)u^{*}
\end{pmatrix}
\begin{pmatrix}
1& a\\
b&-1
\end{pmatrix}
\begin{pmatrix}
1& 1\\
t(w)u&t(w)u
\end{pmatrix}
\begin{pmatrix}
e& 0\\
0&0
\end{pmatrix}
 \\\nonumber
 &=&
 \begin{pmatrix}
e& 0\\
0&0
\end{pmatrix}
\begin{pmatrix}
1+t(w)(u^{*}b+au)-t(w)^{2}& 1+t(w)(u^{*}b+au)-t(w)^{2}\\
1+t(w)(u^{*}b+au)-t(w)^{2}&1+t(w)(u^{*}b+au)-t(w)^{2}
\end{pmatrix}
\begin{pmatrix}
e& 0\\
0&0
\end{pmatrix} \\\nonumber
  &=& 0.
\end{eqnarray}
\end{proof}

\begin{Lemma}\label{alphabeta0}
Let $\M$ be a finite von Neumann algebra and $a, b, c, d\in \M$. If $c$ and $d$ are positive invertible operators, then there is a nonzero operator $y\in M_{2}(\M)$ such that
\[y
\begin{pmatrix}
c& a\\
b&-d
\end{pmatrix}
y^{*}=0.
\]
\end{Lemma}
\begin{proof}
Note that
\[\begin{pmatrix}
c^{-\frac{1}{2}}& 0\\
0&d^{-\frac{1}{2}}
\end{pmatrix}
\begin{pmatrix}
c& a\\
b&-d
\end{pmatrix}
\begin{pmatrix}
c^{-\frac{1}{2}}& 0\\
0&d^{-\frac{1}{2}}
\end{pmatrix}=
\begin{pmatrix}
1& c^{-\frac{1}{2}}ad^{-\frac{1}{2}}\\
d^{-\frac{1}{2}}bc^{-\frac{1}{2}}&-1
\end{pmatrix}.
\]
By lemma \ref{Y1-1Y^{*}}, there is a nonzero operator $\widetilde{y}\in M_{2}(\M)$ such that
\[\widetilde{y}
\begin{pmatrix}
1& c^{-\frac{1}{2}}ad^{-\frac{1}{2}}\\
d^{-\frac{1}{2}}bc^{-\frac{1}{2}}&-1
\end{pmatrix}
\widetilde{y}^{*}
=0.\]
Then \[
\widetilde{y}
\begin{pmatrix}
c^{-\frac{1}{2}}& 0\\
0&d^{-\frac{1}{2}}
\end{pmatrix}
\begin{pmatrix}
c& a\\
b&-d
\end{pmatrix}
\begin{pmatrix}
c^{-\frac{1}{2}}& 0\\
0&d^{-\frac{1}{2}}
\end{pmatrix}
\widetilde{y}^{*}
=0.
\]
Let $y=\widetilde{y}\begin{pmatrix}
c^{-\frac{1}{2}}& 0\\
0&d^{-\frac{1}{2}}
\end{pmatrix}$.
Then $y\neq0$ and
\[y
\begin{pmatrix}
c& a\\
b&-d
\end{pmatrix}
y^{*}=0.
\]
\end{proof}

\begin{Corollary}\label{C-D general}
Let $\M$ be a finite von Neumann algebra and $a, b\in \M$. If $c$ and $d$ are both positive operators, then there is a nonzero operator $y\in M_{2}(\M)$ such that
\[y
\begin{pmatrix}
c& a\\
b&-d
\end{pmatrix}
y^{*}=0.
\]
\end{Corollary}
\begin{proof}
If ker$c\neq \{0\}$. Let $e$ be the projection onto ker$c$. Then $e\neq 0$. Let $y=\begin{pmatrix} e&0\\
0&0
\end{pmatrix}$. Then $y
\begin{pmatrix}
c& a\\
b&-d
\end{pmatrix}
y^{*}=0.$ If ker$d\neq \{0\}$. Let $f$ be the projection onto ker$d$. Then $f\neq 0$. Let $y=\begin{pmatrix} 0&0\\
0&f
\end{pmatrix}$. Then $y
\begin{pmatrix}
c& a\\
b&-d
\end{pmatrix}
y^{*}=0.$ Therefore, we may assume that both $c$ and $d$ are injective with dense range. Note that
\[\begin{pmatrix}
\frac{c^{-\frac{1}{2}}}{\sqrt{1+c^{-1}}}& 0\\
0&\frac{d^{-\frac{1}{2}}}{\sqrt{1+d^{-1}}}
\end{pmatrix}
\begin{pmatrix}
c& a\\
b&-d
\end{pmatrix}
\begin{pmatrix}
\frac{c^{-\frac{1}{2}}}{\sqrt{1+c^{-1}}}& 0\\
0&\frac{d^{-\frac{1}{2}}}{\sqrt{1+d^{-1}}}
\end{pmatrix}=
\begin{pmatrix}
\frac{1}{\sqrt{1+c^{-1}}}& *\\
*&-\frac{1}{\sqrt{1+d^{-1}}}
\end{pmatrix}.
\]
Since $\frac{1}{\sqrt{1+c^{-1}}}$ and $\frac{1}{\sqrt{1+d^{-1}}}$ are invertible, by lemma \ref{alphabeta0}, there is a nonzero operator $x\in M_{2}(\M)$ such that
\[x
\begin{pmatrix}
\frac{1}{\sqrt{1+c^{-1}}}& *\\
*&-\frac{1}{\sqrt{1+d^{-1}}}
\end{pmatrix}
x^{*}=0.
\]
Let \[y=x\begin{pmatrix}
\frac{c^{-\frac{1}{2}}}{\sqrt{1+c^{-1}}}& 0\\
0&\frac{d^{-\frac{1}{2}}}{\sqrt{1+d^{-1}}}
\end{pmatrix}.\] Since $\frac{c^{-\frac{1}{2}}}{\sqrt{1+c^{-1}}}$ and $\frac{d^{-\frac{1}{2}}}{\sqrt{1+d^{-1}}}$ are invertible, it follows that $y\neq0$. Therefore, there is a nonzero operator $y\in M_{2}(\M)$ such that
\[y
\begin{pmatrix}
c& a\\
b&-d
\end{pmatrix}
y^{*}=0.
\]
\end{proof}

\begin{Lemma}\label{finite exe0}
Let $(\M,\tau)$ be a tracial von Neumann algebra and $x\in \M$. If $E_{\mathcal{Z}(\M)}(x)=0$, then there is a nonzero projection $e\in \M$ such that $exe=0$.
\end{Lemma}
\begin{proof}
Write $x=a+ib$, where $a, b\in \M$ are self-adjoint operators. Since $E_{\mathcal{Z}(\M)}(x)=0$, we have \[E_{\mathcal{Z}(\M)}(a)=E_{\mathcal{Z}(\M)}(b)=0.\] By Corollary \ref{8}, there is a sequence of mutually orthogonal nonzero projections $\{e_{n}\}_{n=1}^\infty$ in $\M$ such that $\sum_{n=1}^\infty e_{n}=I$ and $e_{n}ae_{n}=0$ for all $n$.
Now we have
\[
a=\begin{blockarray}{cccccccc}
       & e_{1}   &e_{2}    &\cdots   &e_{n}   &\cdots \\
\begin{block}{c(ccccccc)}
 e_{1} & 0       & a_{12}  &\cdots   &a_{1n}  & \cdots\\
 e_{2} & a_{21}  & 0       &\cdots   &a_{2n}  & \cdots\\
\vdots & \vdots  &\vdots   &\ddots   &\vdots  &\vdots \\
 e_{n} & a_{n1}  &a_{n2}   &\cdots   &0  &\cdots \\
\vdots & \vdots  &\vdots   &\cdots   &\vdots  &\ddots \\
\end{block}
\end{blockarray}
\]
and
\[
b=\begin{blockarray}{cccccccc}
       & e_{1}     &e_{2}    &\cdots   &e_{n}   &\cdots  \\
\begin{block}{c(ccccccc)}
 e_{1} & e_{1}be_{1}     & b_{12}  &\cdots   &b_{1n}  & \cdots  \\
 e_{2} & b_{21}    & e_{2}be_{2}   &\cdots   &b_{2n}  & \cdots  \\
\vdots & \vdots    &\vdots   &\ddots   &\vdots  &\vdots   \\
 e_{n} & b_{n1}    &b_{n2}   &\cdots   &e_{n}be_{n}   &\cdots   \\
\vdots & \vdots    &\vdots   &\cdots   &\vdots  &\ddots \\
\end{block}
\end{blockarray}\,.
\]
For $n\geq 1$, we may assume that $e_{n}be_{n}\geq0$ or $e_{n}be_{n}\leq0$. Indeed, there are subprojections $e_{n1}, e_{n2}$ of $e_n$ such that $e_{n1}+e_{n2}=e_n$ and $e_{n1}be_{n1}\geq 0$ and $e_{n2}be_{n2}\leq 0$. Meanwhile, $e_{n1}ae_{n1}=e_{n2}ae_{n2}=0$. Relabel $\{e_{n1},e_{n2}\}$ as $\{e_n\}$. Then we have $e_nae_n=0$, $\forall n$, and $e_nbe_n\geq 0$ or $e_nbe_n\leq 0$ for all $n$.
Since $E_{\mathcal{Z}(\M)}(b)=0$, we may assume that $e_{1}be_{1}\geq0$ and $e_{2}be_{2}\leq0$.  Note that if $e_{i}\in \mathcal{Z}(\M)$, then \[E_{\mathcal{Z}(e_i\M e_i)}(e_{i}be_{i})=e_{i}E_{\mathcal{Z}(\M)}(b)e_{i}=0.\] By lemma \ref{EXE=0}, there is a nonzero projection $f\in e_{i}\M e_{i}$ such that \[f e_{i}be_{i}f=0.\]So we have a nonzero projection $f\in \M$ such that $f a f=f b f=0$ and thus $fxf=0$.
Let us assume $e_{i}\notin \mathcal{Z}(\M)$ and $e_{i}be_{i}\neq0$ for any $i$. Denote by $\sum_{i\in I}e_{i}=p$ with $e_{i}be_{i}\geq0$, $\sum_{j\in J}f_{j}=q$ with $f_{j}bf_{j}\leq0$, then $p+q=I$.
We claim that there exist projections $e_{i}$ and $f_{j}$ such that $e_{i}z(f_{j})\neq0$. If not, for any $i,j$, we have $e_{i}z(f_{j})=0$. Form this, we also have
\[e_{i}\left(\bigvee_{j\in J}z(f_{j})\right)=0, \quad \forall i.\]
Then
\[\left(\bigvee_{i\in I}z(e_{i})\right)\left(\bigvee_{j\in J}z(f_{j})\right)=0.\] Therefore,  \[ \left(\sum_{i\in I}e_{i}\right)\left(\bigvee_{j\in J}z(f_{j})\right)=0,\quad \left(\bigvee_{i\in I}z(e_{i})\right)\left(\sum_{j\in J}f_{j}\right)=0.\] Since $\sum_{i\in I}e_i+\sum_{j\in J}f_j=I$, we have \[p=\sum_{i\in I}e_{i}=\bigvee_{i\in I}z(e_{i})\in\mathcal{Z}(\M).\] Similarly, \[q=\sum_{j\in J}f_{j}=\bigvee_{j\in J}z(f_{j})\in\mathcal{Z}(\M).\]
Note that
\[0=E_{\mathcal{Z}(\M)}(b)=E_{\mathcal{Z}(\M)}\left(\sum_{i}e_{i}be_{i}+\sum_{j}f_{j}bf_{j}\right).\]
Since $p\in\mathcal{Z}(\M)$, we have \[0=pE_{\mathcal{Z}(\M)}(b)p=E_{\mathcal{Z}(\M)}\left(\sum_{i}e_{i}be_{i}\right).\]
Since $e_{i}be_{i}\geq0$, then we have $e_{i}be_{i}=0$ for any $i$. In the same way, we have $f_{j}bf_{j}=0$ for any $j$.
This contradicts to the assumption that $e_nbe_n\neq 0$ for each $n$. So we may assume that $e_{1}, e_{2}$ are not in $\mathcal{Z}(\M)$ and $e_{1}z(e_{2})\neq0$. Then $e_{1}$ and $e_{2}$ have nonzero equivalent subprojections. Without loss of generality, we may assume that $e_1 \sim e_2$. Let $\widetilde{e}=e_{1}+ e_{2}$ and $\widetilde{b}=\widetilde{e}b\widetilde{e}$. Note that \[(e_{1}+ e_{2})\M(e_{1}+ e_{2})\cong M_{2}(\mathbb{C})\otimes e_{1}\M e_{1}.\] This means that we can represent operators acting on $\widetilde{e}\H$ as $2\times2$ block matrices. Represent $\widetilde{x}=\widetilde{e}x\widetilde{e}$ as a matrix
\[
\widetilde{x}=
\begin{pmatrix}
e_{1}be_{1}& e_{1}ae_{2}+ie_1 b e_2\\
e_2 a e_1+ie_{2}be_{1}& e_{2}be_{2}
\end{pmatrix}=
\begin{pmatrix}
c& a\\
b&-d
\end{pmatrix},
\]
where $c, d\geq0$. By corollary \ref{C-D general} and lemma \ref{11}, there is a nonzero projection $E\in M_2(\mathbb{C})\otimes e_1\M e_1$ such that $E\widetilde{x}E=0$. Let $0\neq e=E\oplus0\in\M$, then $exe=0$.
\end{proof}

\begin{Corollary}\label{C:2.10}
Let $(\mathcal{M},\tau)$ be a tracial von Neumann algebra and $x\in \M$. If $E_{\mathcal{Z}(\M)}(x)=0$, then there is a sequence of nonzero projections $\{p_n\}_{n=1}^\infty$ in $\M$ such that $\sum_{n=1}^\infty p_n=I$ and $p_nxp_n=0$ for all $n\geq 1$.
\end{Corollary}
\begin{proof}
Let $\E=\{\{p_\alpha\}\}$ be the family of mutually orthogonal nonzero projections such that $p_\alpha xp_\alpha=0$ for all $\alpha$. By lemma \ref{finite exe0}, $\E$ is not empty. Define $\{p_\alpha\}\prec \{q_\beta\}$ if $\{p_\alpha\}\subseteq \{q_\beta\}$. Clearly, $\prec$ is a partial order and if $\{\{p_\alpha\}\}$ is a totally ordered set, then the union is a maximal element of the totally ordered set. By Zorn's lemma, there is a maximal element $\{p_n\}$ in $\E$. We claim that $\sum_n p_n=I$. Otherwise, let $f=I-\sum_{n}p_n$.
Note that
\begin{align*}
 E_{\Z(\M)}(fxf)&=E_{\Z(\M)}(xf)=E_{\Z(\M)}\left(x-x\sum_n p_n\right)\\
               &=E_{\Z(\M)}(x)-\sum_n  E_{Z(\M)}(p_n xp_n)=0.
\end{align*}
Since
\begin{eqnarray}\nonumber
  E_{\Z(\M)}(fxf)=0 &\Leftrightarrow& \langle fxf,y\rangle_{\tau}=0, \, \forall y \in \Z(\M), \\  \nonumber
   &\Leftrightarrow& \langle fxf,fyf\rangle_{\tau}=0, \, \forall y \in \Z(\M), \\  \nonumber
   &\Leftrightarrow& E_{f\Z(\M)f}(fxf)=0,
\end{eqnarray}
then, by lemma \ref{finite exe0}, there is a nonzero projection $p$ in $f\M f$ such that \[pfxfp=pxp=0.\]
Now $\{p_n\}\cup \{p\}\in \E$. This is a contradiction.

\end{proof}

\begin{Lemma}\label{7'}
Let $\M$ be a type ${\rm II}_1$ von Neumann algebra and $p_1, q_1$ be two nonzero equivalent orthogonal projections in $\M$. Let \[
H_i=\begin{pmatrix}
0&\quad d_i\\
d_i^*&\quad 0
\end{pmatrix}\begin{matrix}
p_1\\
q_1
\end{matrix}
, \quad K=\begin{pmatrix}
a &\quad c\\
c^*&\quad -b
\end{pmatrix}\begin{matrix}
p_1\\
q_1
\end{matrix},
\]
where $a,b$ are positive operators and $1\leq i\leq n$. Then there is a nonzero projection $f\leq p_1+q_1$ such that $fH_if=fKf=0$, $1\leq i\leq n$.

\end{Lemma}
\begin{proof}
Since $p_1$ and $q_1$ are equivalent projections, then
\[(p_1+q_1)\M (p_1+q_1)\cong M_2(\mathbb{C})\otimes p_1\M p_1.\]

Assume $n=1$. Consider the operator
\[
K+iH_1=\begin{pmatrix}
a &\quad \ast\\
\ast &\quad -b
\end{pmatrix}\begin{matrix}
p_1\\
q_1
\end{matrix}.
\]
By corollary~\ref{C-D general} and lemma~\ref{11}, there is a nonzero projection $f\leq p_1+q_1$ such that $f(K+iH_1)f=0$. Therefore,  $fH_1f=fKf=0$.

Assume the lemma is true for $n=k$.
By the polar decomposition theorem, $d_1=u|d_1|$, where $u$ is a unitary operator. Thus
\[
\begin{pmatrix}
0& d_1   \\
d_1^*& 0
\end{pmatrix}
=\begin{pmatrix}
u& 0   \\
0& I
\end{pmatrix}
\begin{pmatrix}
0& |d_1|   \\
|d_1|& 0
\end{pmatrix}
\begin{pmatrix}
u^*& 0   \\
0& I
\end{pmatrix}
\]
and
\[
\begin{pmatrix}
u^*& 0   \\
0& I
\end{pmatrix}
\begin{pmatrix}
a& c   \\
c^*& -b
\end{pmatrix}
\begin{pmatrix}
u& 0   \\
0& I
\end{pmatrix}
=
\begin{pmatrix}
u^*au& u^*c   \\
c^*u & -b
\end{pmatrix}.
\]
So we may assume that $d_1\geq 0$. Then there is a maximal abelian von Neumann algebra $\A$ of $p_1\M p_1$ containing $d_1$.
Since $p_1\M p_1$ is a type ${\rm II}_1$ von Neumann algebra, there is a projection $e\in \A$ such that $e$ is not in $\Z(p_1\M p_1)$. Thus there exists a central projection $\tilde{e}$ in $p_1\M p_1$ such that $e\tilde{e}\preceq (p_1-e)\tilde{e}$ and $(p_1-e)(p_1-\tilde{e})\preceq e(p_1-\tilde{e})$. Since $e$ is not in $\Z(p_1\M p_1)$, then we have either $e\tilde{e}\neq 0$ or $(p_1-e)(p_1-\tilde{e})\neq 0$. We may assume that $e\tilde{e}\neq 0$. Let $p=e\tilde{e}\in \A$ and $q=p_1-p$. Then $p\preceq q$. Let $\tilde{q}$ be a subprojection of $q$ such that $p\sim \tilde{q}$ in $p_1\M p_1$. Note that $\tilde{q}\notin \A$ in general. However, we still have $pd_1\tilde{q}=\tilde{q}d_1 p=0$. Thus we can write
\[
H_1=
\begin{pmatrix}
0 & 0&  pd_{1}&0\\
0 & 0&  0&\tilde{q}d_{1}\tilde{q}\\
pd_{1} & 0&  0&0\\
0 & \tilde{q}d_{1}\tilde{q}&  0&0
\end{pmatrix},\quad
H_i=
\begin{pmatrix}
0 & 0&  \ast&\ast\\
0 & 0&  \ast &\ast\\
\ast & \ast&  0&0\\
\ast & \ast &  0&0
\end{pmatrix}
\]
for $2\leq i\leq k+1$ and
\[
K=
\begin{pmatrix}
pap & *&  *&*\\
* & \tilde{q}a\tilde{q}&  *&*\\
* & *&  -pbp&*\\
* & *&  *&-\tilde{q}b\tilde{q}
\end{pmatrix}.
\]
Considering the second and third columns and rows of $H_i$ and $K$, by the mathematical induction,  there is a nonzero projection $f\leq p+\tilde{q}$ such that $fH_if=fKf=0$, $2\leq i\leq k+1$. Note that the second and third columns and rows of $H_1$ is a zero matrix. Therefore, we have $fH_1f=0$ and complete the proof.

\end{proof}

\begin{Lemma}\label{L:center}
Let $(\mathcal{M}, \tau)$ be a type ${\rm II}_1$ von Neumann algebra. For any finite self-adjoint operators $x_{1},x_{2},\ldots,x_{m}\in\mathcal{M}$, $E_{\Z(\mathcal{M})}(x_{j})=0$, $1\leq j\leq m$ and any $\epsilon>0$, there exists a partition of 1 with projections $p_{1},\ldots, p_{n}$ in $\mathcal{M}$ such that
\[\left\|\sum_{i=1}^n p_{i}x_{j}p_{i}\right\|_{2}<\epsilon, \quad 1\leq j\leq m.\]
\end{Lemma}
\begin{proof}
We only need to show that there is a sequence of projections $p_{1},\ldots, p_{n},\ldots$ such that $\sum p_n=I$ and $p_nx_jp_n=0$ for $1\leq j\leq m$. The case $m=1$ follows from corollary~\ref{8} and the case $m=2$ follows from corollary~\ref{C:2.10}. We prove the case $m=3$. First, we show that there is a nonzero projection $p\in \M$ such that $px_ip=0$ for $1\leq i\leq 3$. We can write
\[
x_{i}=
\bordermatrix	
{      &q_1    &q_2    &\cdots &q_n    &\cdots \cr
q_1    &0      &\ast   &\cdots &\ast   &\cdots \cr
q_2    &\ast   & 0     &\cdots &\ast   &\cdots \cr
~\vdots&\vdots &\vdots &\ddots &\vdots &\vdots \cr
q_n    &\ast   &\ast   &\cdots &0      &\cdots \cr
~\vdots&\vdots &\vdots &\vdots &\vdots &\ddots \cr},\quad i=1,2,\]
\[
x_{3}=
\bordermatrix	
{       &q_1         &q_2         &\cdots &q_n         &\cdots \cr
q_1     &q_1x_{3}q_1 &\ast        &\cdots &\ast        &\cdots \cr
q_2     &\ast        &q_2x_{3}q_2 &\cdots &\ast        &\cdots \cr
~\vdots &\vdots      &\vdots      &\ddots &\vdots      &\vdots \cr
q_n     &\ast        &\ast        &\cdots &q_nx_{3}q_n &\cdots \cr
~\vdots &\vdots      &\vdots      &\vdots &\vdots      &\ddots \cr}
\]
for a decomposition $I=q_1+\cdots+q_n+\cdots$. Since $E_{\Z(\M)}(x_{3})=0$, we may assume that $q_k x_3q_k$ is positive or negative for all $k\geq 1$. Let us assume that $q_1 x_{3}q_1\geq 0$ and $q_2 x_{3}q_2\leq 0$. We may assume that $q_1$ and $q_2$ are not in $\Z(\M)$. Then by cutting $q_1$ and $q_2$ suitably, we may assume that $q_1$ and $q_2$ are equivalent projections in $\M$. Now by lemma~\ref{7'}, there is a nonzero projection $p\in \M$ such that $px_1p=px_2p=px_3p=0$.

Let $\E=\{\{e_\alpha\}\}$ be the family of mutually orthogonal nonzero projections such that \[e_\alpha x_1e_\alpha=e_\alpha x_2 e_\alpha=e_\alpha x_3 e_\alpha=0\] for all $\alpha$. By the above statetment, then $\E$ is not empty. Define $\{e_\alpha\}\prec \{f_\beta\}$ if $\{e_\alpha\}\subseteq \{f_\beta\}$. Clearly, $\prec$ is a partial order and if $\{\{e_\alpha\}\}$ is a totally ordered set, then the union is a maximal element of the totally ordered set. By Zorn's lemma, there is a maximal element $\{e_n\}_{n=1}^\infty$ in $\E$. We claim that $\sum_{n=1}^\infty e_n=I$. Otherwise, let $f=I-\sum_{n=1}^\infty e_n$. Then \[E_{\Z(f\M f)}(fx_if)=0,\, 1\leq i\leq 3.\] Thus there is a nonzero subprojection $e$ of $f$ such that $ex_1e=ex_2e=ex_3e=0$. Then $\{e_n\}\cup \{e\}\in \E$. It is a contradiction.

The general case can be proved by lemma~\ref{7'} and mathematical induction.

\end{proof}

\begin{Remark}
Let $(\mathcal{M}, \tau)$ be a tracial von Neumann algebra. For any finite self-adjiont operators $x_{1},x_{2},\ldots,x_{m}\in\mathcal{M}$ with $m \geq 3$, $E_{\Z(\mathcal{M})}(x_{j})=0$, $1\leq j\leq m$, there may not exists a nonzero projection $p$ in $\mathcal{M}$ such that
\[ p x_{j}p =0, \quad 1\leq j\leq m.\]
For example, let $\mathcal{M}=M_2(\mathbb{C})$ and
\[x_1=
\begin{pmatrix}
0& 1   \\
1& 0
\end{pmatrix},\,\,\,
x_2=
\begin{pmatrix}
0& -i   \\
i& 0
\end{pmatrix},\,\,\,
x_3=
\begin{pmatrix}
1& 0   \\
0& -1
\end{pmatrix}.
\]
Then $\tau(x_1)=\tau(x_2)=\tau(x_3)=0$. But there does not exist nonzero projection $p$ such that $px_1p=px_2p=px_3p=0$. For atomic tracial von Neumann algebra, the above statement is not true. Then we have the following natural question.

\end{Remark}

\textbf{Question:}
Let $(\mathcal{M}, \tau)$ be a diffuse tracial von Neumann algebra. For any finite set of operators $x_{1},x_{2},\ldots,x_{m}\in\mathcal{M}$, $E_{\Z(\mathcal{M})}(x_{j})=0$, $1\leq j\leq m$ and any $\epsilon>0$, does there exists a sequence of projections $\{p_n\}_{n=1}^{\infty}$ in $\M$ such that $\sum_{n=1}^{\infty}p_n=I$ and
\[ p_n x_j p_n=0 \quad\text{for}\quad n \geq 1, \,1\leq j \leq m \,?  \]

If the above statement holds true, according to the proof of the Theorem \ref{T:local quantization principle} below, we know that Theorem~\ref{T:local quantization principle} also holds for diffuse tracial von Neumann algebra.
\vskip 0.5cm

The following lemma is due to Sorin Popa (see Lemma A.1.1 of~\cite{popa 94}).
\begin{Lemma}{\label{7}}
Let $\mathcal{N}\subseteq \mathcal{M}$ be finite von Neumann algebras with a normal finite faithful trace $\tau$. Let $X\subset\mathcal{M}$ be a finite set of elements such that $E_{\mathcal{N}\vee(\mathcal{N}'\cap \mathcal{M})}(x)=0$ for all $x\in X$. Given any $ \epsilon> 0$, there exists a partition of the identity with projections $\{p_i\}_{i=1}^n$ in $\mathcal{N}$ such that $\|\sum_{i=1}^np_ixp_i\|_{2}<\epsilon$ for all $x\in X $.
\end{Lemma}

\begin{Lemma}\label{6}
Let $(\mathcal{M},\tau)$ be a type ${\rm II}_1$ von Neumann algebra  and $\mathcal{N}\subseteq \mathcal{M}$ a type ${\rm II}_1$ von Neumann subalgebra. For any $x\in \mathcal{M}$ and any $ \epsilon> 0$, $E_{\mathcal{N}'\cap \mathcal{M}}(x)=0$, then there exists a partition of 1 with projections $p_{1}, \ldots, p_{n}$ in $\mathcal{N}$ such that \[\left\|\sum_{i=1}^n p_{i}xp_{i}\right\|_{2}<\epsilon.\]
\end{Lemma}
\begin{proof}
Let $y=E_{\mathcal{N}\vee(\mathcal{N}'\cap \mathcal{M})}(x)$. Then $y\in\mathcal{N}\vee(\mathcal{N}'\cap \mathcal{M})$ and \[E_{\mathcal{N}'\cap \mathcal{M}}(y)=E_{\mathcal{N}'\cap \mathcal{M}}(E_{\mathcal{N}\vee(\mathcal{N}'\cap \mathcal{M})}(x))=E_{\mathcal{N}'\cap \mathcal{M}}(x)=0.\] Note that $x=x-y+y$ and $E_{\mathcal{N}\vee(\mathcal{N}'\cap \mathcal{M})}(x-y)=0$. For each $y$ as above there exist  $x_{1},x_{2},\ldots,x_{N}\in \mathcal{N}$, $x_{1}',x_{2}',\ldots,x_{N}'\in \mathcal{N}'\cap \mathcal{M}$ such that \[\left\|y-\sum_{i=1}^{N}x_{i}x_{i}'\right\|_{2}<\frac{\epsilon}{2}.\] Then
\[\left\|E_{\mathcal{N}'\cap \mathcal{M}}(y)-\sum_{i=1}^{N}E_{\mathcal{N}'\cap \mathcal{M}}(x_{i})x_{i}'\right\|_{2}<\frac{\epsilon}{2}.\]
Note that $E_{\N'\cap \M}(y)=0$ and
\[E_{\mathcal{N}'\cap \mathcal{M}}(x_{i})=E_{\mathcal{N}'\cap \mathcal{M}}(E_{\mathcal{N}}(x_{i}))=E_{\Z(\mathcal{N})}(x_{i}).\]
Then \[\left\|\sum_{i=1}^{N}E_{\Z(\mathcal{N})}(x_{i})x_{i}'\right\|_{2}<\frac{\epsilon}{2}\] and \[\left\|y-\sum_{i=1}^{N}(x_{i}-E_{\Z(\mathcal{N})}(x_{i}))x_{i}'\right\|_{2}\leq\left\|y-\sum_{i=1}^{N}x_{i}x_{i}'\right\|_{2}+
\left\|\sum_{i=1}^{N}E_{\Z(\mathcal{N})}(x_{i})x_{i}'\right\|_{2}<\epsilon.\] Since \[E_{\Z(\mathcal{N})}(x_{i}-E_{\Z(\mathcal{N})}(x_{i}))=0,\] by lemma \ref{L:center}, there exists a partition of 1 with projections $q_{1}, \ldots, q_{m}$ in $\mathcal{N}$ such that
\[\left\|\sum_{j=1}^m q_{j}(x_{i}-E_{\Z(\mathcal{N})}(x_{i}))q_{j}\right\|_{2}<\frac{\epsilon}{N\|x_{i}'\|}.\] Then \[\left\|\sum_{j=1}^m q_{j}\left(\sum_{i=1}^{N}(x_{i}-E_{\Z(\mathcal{N})}(x_{i}))x_{i}'\right)q_{j}\right\|_{2}\leq\sum_{i=1}^{N}\left\|\sum_{j=1}^m q_{j}(x_{i}-E_{\Z(\mathcal{N})}(x_{i}))q_{j}\right\|_{2}\|x_{i}'\|<\epsilon\] and
\begin{eqnarray}\nonumber
   \left\|\sum_{j=1}^m q_{j}yq_{j}\right\|_{2}&\leq& \left\|\sum_{j=1}^m q_{j}yq_{j}-\sum_{j=1}^m q_{j}\left(\sum_{i=1}^{N}(x_{i}-E_{\Z(\mathcal{N})}(x_{i}))x_{i}'\right)q_{j}\right\|_{2} \\\nonumber
   && +\left\|\sum_{j=1}^m q_{j}\left(\sum_{i=1}^{N}(x_{i}-E_{\Z(\mathcal{N})}(x_{i}))x_{i}'\right)q_{j}\right\|_{2} \\\nonumber
   &<&  2\epsilon.
\end{eqnarray}
Note that \[(q_{j}\mathcal{N}q_{j})'\cap q_{j}\mathcal{M}q_{j}=q_{j}(\mathcal{N}'\cap\mathcal{M})q_{j}\] and \[q_{j}\mathcal{N}q_{j}\vee((q_{j}\mathcal{N}q_{j})'\cap q_{j}\mathcal{M}q_{j})=q_j(\mathcal{N}\vee(\mathcal{N}'\cap \mathcal{M}))q_j, \quad 1\leq j\leq m.\] Since $q_{j}\in\mathcal{N}\subseteq\mathcal{N}\vee(\mathcal{N}'\cap \mathcal{M})$, then
\begin{eqnarray}\nonumber
  E_{q_{j}\mathcal{N}q_{j}\vee((q_{j}\mathcal{N}q_{j})'\cap q_{j}\mathcal{M}q_{j})}(q_j(x-y)q_j) &=&  E_{q_j(\mathcal{N}\vee(\mathcal{N}'\cap \mathcal{M}))q_j}(q_j(x-y)q_j)\\\nonumber
   &=& q_j E_{\mathcal{N}\vee(\mathcal{N}'\cap \mathcal{M})}(x-y)q_j \\\nonumber
   &=&0,\quad 1\leq j\leq m.
\end{eqnarray}
By lemma \ref{7}, there exists a partition of $q_{j}$ with projections $p_{i_1}, \ldots, p_{i_{m_j}}$ in $\mathcal{N}$ such that \[\left\|\sum_{k=1}^{m_j}p_{i_k}(q_j(x-y)q_j)p_{i_k}\right\|_{2}=\left\|\sum_{k=1}^{m_j}p_{i_k}(x-y)p_{i_k}\right\|_{2}<\frac{\epsilon}{m}. \] Note that \[\left\|\sum_{j=1}^m\sum_{k=1}^{m_j} p_{i_k}(q_{j}yq_{j})p_{i_k}\right\|_{2}\leq\left\|\sum_{j=1}^mq_{j}yq_{j}\right\|_{2}<2\epsilon.\]
 Relabel $\{p_{i_k}\}$ as $\{p_r\}$, $1\leq r\leq n$. Then $\sum_{i=1}^np_i=I$ and \[\left\|\sum_{i=1}^{n}p_ixp_i\right\|_{2}=\left\|\sum_{i=1}^{n}p_i(x-y)p_i\right\|_{2}+\left\|\sum_{i=1}^{n}p_iyp_i\right\|_{2}<3\epsilon.\]
\end{proof}

\noindent{{\bf Proof of Theorem~\ref{T:local quantization principle}}}: Clearly we may assume that $E_{\N'\cap \M}(x_i)=0$ for $1\leq i\leq m$. Let $y_i=E_{\mathcal{N}\vee(\mathcal{N}'\cap \mathcal{M})}(x_i)$, $1\leq i\leq m$. Then $y_i\in\mathcal{N}\vee(\mathcal{N}'\cap \mathcal{M})$ and \[E_{\mathcal{N}'\cap \mathcal{M}}(y_i)=E_{\mathcal{N}'\cap \mathcal{M}}(E_{\mathcal{N}\vee(\mathcal{N}'\cap \mathcal{M})}(x_i))=E_{\mathcal{N}'\cap \mathcal{M}}(x_i)=0,\quad 1\leq i\leq m.\] Note that $x_i=x_i-y_i+y_i$ and $E_{\mathcal{N}\vee(\mathcal{N}'\cap \mathcal{M})}(x_i-y_i)=0$. For each $y_i$ as above there exist  $x_{i1},x_{i2},\ldots,x_{iN}\in \mathcal{N}$, $x_{i1}',x_{i2}',\ldots,x_{iN}'\in \mathcal{N}'\cap \mathcal{M}$ such that \[\left\|y_i-\sum_{k=1}^{N}x_{ik}x_{ik}'\right\|_{2}<\frac{\epsilon}{2}, \quad 1\leq i\leq m.\] Then
\[\left\|E_{\mathcal{N}'\cap \mathcal{M}}(y_i)-\sum_{k=1}^{N}E_{\mathcal{N}'\cap \mathcal{M}}(x_{ik})x_{ik}'\right\|_{2}<\frac{\epsilon}{2}.\]
Note that $E_{\mathcal{N}'\cap \mathcal{M}}(y_i)=0$ and
\[E_{\mathcal{N}'\cap \mathcal{M}}(x_{ik})=E_{\mathcal{N}'\cap \mathcal{M}}(E_{\mathcal{N}}(x_{ik}))=E_{\Z(\mathcal{N})}(x_{ik}).\]
Then \[\left\|\sum_{k=1}^{N}E_{\Z(\mathcal{N})}(x_{ik})x_{ik}'\right\|_{2}<\frac{\epsilon}{2}\] and \[\left\|y_i-\sum_{k=1}^{N}(x_{ik}-E_{\Z(\mathcal{N})}(x_{ik}))x_{ik}'\right\|_{2}\leq\left\|y_i-\sum_{k=1}^{N}x_{ik}x_{ik}'\right\|_{2}+
\left\|\sum_{k=1}^{N}E_{\Z(\mathcal{N})}(x_{ik})x_{ik}'\right\|_{2}<\epsilon,\] for $1\leq i\leq m$. Since $E_{\Z(\mathcal{N})}(x_{ik}-E_{\Z(\mathcal{N})}(x_{ik}))=0$, by lemma \ref{L:center}, there exists a partition of 1 with projections $q_{1}, \ldots, q_{l}$ in $\mathcal{N}$ such that
\[\left\|\sum_{j=1}^l q_{j}(x_{ik}-E_{\Z(\mathcal{N})}(x_{ik}))q_{j}\right\|_{2}<\frac{\epsilon}{N\|x_{ik}'\|}.\] Then \[\left\|\sum_{j=1}^l q_{j}\sum_{k=1}^{N}(x_{ik}-E_{\Z(\mathcal{N})}(x_{ik}))x_{ik}'q_{j}\right\|_{2}\leq\sum_{k=1}^{N}\left\|\sum_{j=1}^l q_{j}(x_{ik}-E_{\Z(\mathcal{N})}(x_{ik}))q_{j}\right\|_{2}\|x_{ik}'\|<\epsilon.\] and
\begin{eqnarray}\nonumber
   \left\|\sum_{j=1}^l q_{j}y_iq_{j}\right\|_{2}&\leq& \left\|\sum_{j=1}^l q_{j}y_iq_{j}-\sum_{j=1}^l q_{j}\left(\sum_{k=1}^{N}(x_{ik}-E_{\Z(\mathcal{N})}(x_{ik}))x_{ik}'\right)q_{j}\right\|_{2} \\\nonumber
   && +\left\|\sum_{j=1}^l q_{j}\left(\sum_{k=1}^{N}(x_{ik}-E_{\Z(\mathcal{N})}(x_{ik}))x_{ik}'\right)q_{j}\right\|_{2} \\\nonumber
   &<&  2\epsilon, \quad 1\leq j\leq l.
\end{eqnarray}
Note that \[(q_{j}\mathcal{N}q_{j})'\cap q_{j}\mathcal{M}q_{j}=q_{j}(\mathcal{N}'\cap\mathcal{M})q_{j}\] and \[q_{j}\mathcal{N}q_{j}\vee((q_{j}\mathcal{N}q_{j})'\cap q_{j}\mathcal{M}q_{j})=q_j(\mathcal{N}\vee(\mathcal{N}'\cap \mathcal{M}))q_j, \quad 1\leq j\leq l.\] Since $q_{j}\in\mathcal{N}\subseteq\mathcal{N}\vee(\mathcal{N}'\cap \mathcal{M})$, then
\begin{eqnarray}\nonumber
  E_{q_{j}\mathcal{N}q_{j}\vee((q_{j}\mathcal{N}q_{j})'\cap q_{j}\mathcal{M}q_{j})}(q_j(x_i-y_i)q_j) &=&  E_{q_j(\mathcal{N}\vee(\mathcal{N}'\cap \mathcal{M}))q_j}(q_j(x_i-y_i)q_j)\\\nonumber
   &=& q_j E_{\mathcal{N}\vee(\mathcal{N}'\cap \mathcal{M})}(x_i-y_i)q_j \\\nonumber
   &=&0,\quad 1\leq j\leq l
\end{eqnarray}
By lemma \ref{7}, there exists a partition of $q_{j}$ with projections $p_{j_1}, \ldots, p_{j_{m_j}}$ in $\mathcal{N}$ such that \[\left\|\sum_{r=1}^{m_j}p_{j_r}(q_j(x_i-y_i)q_j)p_{j_r}\right\|_{2}=\left\|\sum_{r=1}^{m_j}p_{j_r}(x_i-y_i)p_{j_r}\right\|_{2}<\frac{\epsilon}{l},\quad 1\leq j\leq l. \] Note that \[\left\|\sum_{j=1}^l\sum_{r=1}^{m_j} p_{j_r}(q_{j}y_iq_{j})p_{j_r}\right\|_{2}\leq\left\|\sum_{j=1}^{l}q_jy_iq_{j}\right\|_{2}<2\epsilon,\quad 1\leq i\leq m.\]
Relabel $p_{j_r}$ as $p_k$, $1\leq k\leq n$. Then \[\left\|\sum_{k=1}^{n}p_kx_ip_k\right\|_{2}=\left\|\sum_{k=1}^{n}p_k(x_i-y_i)p_k\right\|_{2}+\left\|\sum_{k=1}^{n}p_ky_ip_k\right\|_{2}<3\epsilon,\quad 1\leq i\leq m.\]
This ends the proof.

\qed
\vskip 1cm

\section{Proof of Theorem~\ref{index 2 finite proj in N}}
In this section, let $\N\subseteq \M$ be an inclusion of type $\rm II_{1}$ factors such that $[\M:\N]=2$. To prove Theorem~\ref{index 2 finite proj in N}, we only need to prove the following theorem.
\begin{Theorem}\label{index 2 finite proj in N cut to 0 many trace}
Let $\N\subseteq \M$ be an inclusion of type $\rm II_{1}$ factors and $[\M:\N]=2$. For any $x_{1},\ldots, x_{m}\in \M$ with $\tau(x_j)=0,~j=1,\ldots,m$, there is a family of finitely many mutually orthogonal nonzero projections $\{e_i\}_{i=1} ^{N}$ in $\mathcal{N}$ such that $\sum_{i=1}^{N}e_i=I$ and $e_ix_je_i=0$ for $1\leq j\leq m$ and $1\leq i\leq N$.
\end{Theorem}

To prove Theorem \ref{index 2 finite proj in N cut to 0 many trace}, we need the following lemmas.
\begin{Lemma}{\text{{\rm (}Goldman's theorem in \cite{Goldman}{\rm )}}}
Let $\N$ be a factor of the $\rm {II}_1$ factor $\M$ with $[\M:\N]=2$. Then $\M$ decomposes as the crossed product of $\N$ by an outer action of $\mathbb{Z}_2$.
\end{Lemma}

\begin{Lemma}\label{outer automorphism}{\rm (}Theorem 1.2.1 in \cite{Alain Connes}{\rm )}
Let $\M$ be a countably decomposable von Neumann algbera and $\theta\in \text{Aut}(\M)$. Then $\theta$ is properly outer if and only if for any nonzero projection $e\in \M$ and any $\epsilon >0$ there exists a nonzero projection $f\leq e$ such that $\|f\theta(f)\|\leq\epsilon$.
\end{Lemma}

\begin{Lemma}\label{WFY3}{\rm (}Theorem 1.1 in \cite{WFY3}{\rm )}
Let $(\mathcal{M},\tau)$ be a type ${\rm II}_1$ factor and $x_1,\ldots, x_n \in \mathcal{M}$. Then there is a family of finitely many mutually orthogonal nonzero projections $\{e_i \}_{i=1} ^{N}$ in $\mathcal{M}$ such that $\sum_{i=1} ^{N}e_i =I$ and $e_ix_je_i=\tau (x_j)e_i$ for $1\leq i\leq N$ and $1\leq j\leq n$.
Equivalently, there is a unitary operator $u\in \M$ such that \[\frac{1}{N} \sum_{i=0}^{N-1}u^{\ast i} x_j u^i=\tau(x_j)I\] for $1\leq j\leq n$.
\end{Lemma}

The following result plays a key role in the proof of the main theorem of this section.

\begin{Lemma}\label{two projection}
Let $\M=\N\rtimes \mathbb{Z}_{2}$ and $L_{g}\in \M$ be the unitary operator corresponding to the generator of $\mathbb{Z}_{2}$. Then there are nonzero projections $e$ and $f$ in $\N$, $e+f=I$, such that $eL_{g}e=fL_{g}f=0$.
\end{Lemma}

To prove Lemma \ref{two projection}, we need the following observations.

Let $\M=\N\rtimes \mathbb{Z}_{2}$. The action of $\mathbb{Z}_{2}$ on $\N$ is given by an automorphism $\sigma_{g}:\N\rightarrow \N$ such that  $\sigma_{g}^{2}=id$. For any nonzero operator $x\in \N$, we have \[\sigma_{g}(x+\sigma_{g}(x))=x+\sigma_{g}(x)\] and \[\sigma_{g}(x-\sigma_{g}(x))=-(x-\sigma_{g}(x)).\]
Then we decompose $x$ as:
\[x=\frac{x+\sigma_{g}(x)}{2}+\frac{x-\sigma_{g}(x)}{2}.\]
Since \[[\N:\N_1]=2,\] where $\N_1$ denotes the fixed point algebra of $\sigma_{g}$, there exists a nonzero operator $x=x^{*}\in \N$ such that $\sigma_{g}(x)=-x$. Let $x=vh$ be the polar decomposition of $x$. Note that \[x=x^{*}=hv^{*}=v^{*}vhv^{*}.\] Then $v=v^{*}$ and $hv=vh$. The automorphism $\sigma_g$ acts on $x$ satisfying \[\sigma_g(x)=\sigma_g(v)\sigma_g(h)=-vh.\] The uniqueness of the polar decomposition implies that \[\sigma_g(v)=-v ~~\text{and}~~ \sigma_g(h)=h.\] Let $p=v^{2}$. Then \[\sigma_{g}(p)=\sigma_{g}(v^{2})=v^{2}=p.\]

\noindent\textbf{Claim:} there exists a unitary operator $u\in \N$ with  $u=u^{*}$ such that \[\sigma_{g}(u)=-u.\]
In fact, let
\[\mathcal{J}=\{v\in \N~|~v=v^{*}~ \text{be a partial isometry and}~ \sigma_{g}(v)=-v\}.\]
It is easy to see that $\mathcal{J}$ is non-empty. Define a relation on $\mathcal{J}$: \[v_{1}\leq v_{2}\Leftrightarrow v_{1}^{2}=e_{1}\leq e_{2}=v_{2}^{2}.\] The relation ``$\leq$'' is a partial order on $\mathcal{J}$. Given any total order subset $\{v_{\alpha}\}$ of $\mathcal{J}$, let \[v=\text{SOT-lim}~v_{\alpha}.\]
Since $v_{\alpha}=v_{\alpha}^{*}$ for any $\alpha$, then \[v^{*}=\text{SOT-lim}~v_{\alpha}^{*}=\text{SOT-lim}~v_{\alpha}=v.\]
For any $\alpha$, $\sigma_{g}(v_{\alpha})=-v_{\alpha}$. Since $\sigma_{g}$ is SOT-continuous, we have
\[\sigma_{g}(v)=\sigma_{g}(\text{SOT-lim}~v_{\alpha})=
\text{SOT-lim}~\sigma_{g}(v_{\alpha})=\text{SOT-lim}(-v_{\alpha})=-v.\]
Thus $v\in\mathcal{J}$. By Zorn's lemma, there is a maximal element $u$ in $\mathcal{J}$. Next, we prove that $u$ is a unitary operator, that is, $u^{2}=I$. Now, we have \[\sigma_{g}(u)=-u ~~\text{and}~~ \sigma_{g}(u^{2})=u^{2}.\] If $u$ is not a unitary operator, let $e =I-u^{2}\neq0$, then \[u^{2}=I-e\] and \[\sigma_{g}(I-e)=\sigma_{g}(u^{2})=u^{2}=I-e.\] By lemma \ref{outer automorphism}, we have $\sigma_g|_{e\N e}\neq id$. Then there is a nonzero partial isometry $v\in e\N e$ with $v=v^*$ such that $\sigma_{g}(v)=-v$. Let
\[w=u+v.\]
Since

(1) $w^{*}=(u+v)^{*}=u^{*}+v^{*}=u+v=w$;

(2) Since $w^{*}w=(u+v)^{*}(u+v)=u^{2}+uv+vu+v^{2}=1-e+v^2$ is a projection, then $w$ is a partial isometry;

(3) $\sigma_{g}(w)=\sigma_{g}(u+v)=\sigma_{g}(u)+\sigma_{g}(v)=-(u+v)$.\\

Therefore, $w\in \mathcal{J}$, contradicting with the maximality of $u$. Thus we proved the \textbf{Claim}.

\qed

\noindent\textbf{The proof of Lemma \ref{two projection}:}
By the previous statement, let $u=2e-1$ be a unitary operator such that $\sigma_{g}(u)=-u$. Note that \[\tau(u)=\tau(\sigma_{g}(u))=-\tau(u).\] Then $\tau(u)=0$ and $\tau(e)=\frac{1}{2}$. Since \[\sigma_{g}(2e-1)=\sigma_{g}(u)=-u=-(2e-1),\] then \[\sigma_{g}(e)=1-e \quad\text{and}\quad \sigma_{g}(1-e)=e.\] Therefore
\[eL_{g}e=e\sigma_{g}(e)L_{g}=e(1-e)L_{g}=0\] and
\[(1-e)L_{g}(1-e)=(1-e)\sigma_{g}((1-e))L_{g}=(1-e)e L_{g}=0.\]
Let $f=1-e$. This completes the proof.

\qed

\begin{Lemma}\label{sigma b=-b}
Let $\M=\N\rtimes \mathbb{Z}_{2}$ and $b=b^*\in \N$ with $\sigma_g(b)=-b$. Then there are nonzero projections $e$ and $f$ in $\N$ with $e+f=I$ and $\sigma_g(e)=f$, such that $ebL_{g}e=fbL_{g}f=0$.
\end{Lemma}
\begin{proof}
Let $b=vh$ be the polar decomposition of $b$. Note that \[b=b^{*}=hv^{*}=v^{*}vhv^{*}.\] Then $v=v^{*}$ and $hv=vh$. Since \[\sigma_g(b)=\sigma_g(v)\sigma_g(h)=-vh,\]
and the uniqueness of the polar decomposition, we have \[\sigma_g(v)=-v ~~\text{and}~~ \sigma_g(h)=h.\] Let $p=v^{2}$. Then \[\sigma_{g}(p)=\sigma_{g}(v^{2})=v^{2}=p.\] Write $v=2e_1-p$, where $0\leq e_1\leq p$. Since $\sigma_g(v)=-v$, we have \[\sigma_g(e_1)=p-e_1.\]
Since $e_{1}=\frac{1}{2}(v+p)=\frac{1}{2}(v+v^{2})$, we have
\[e_{1}b=\frac{1}{2}(v+v^{2})vh=\frac{1}{2}(v+v^{2})hv=\frac{1}{2}vh(v+v^{2})=be_{1}.\]
Note that $(p-e_1)e_1=e_1(p-e_1)=0$. Thus
\[e_1bL_ge_1=e_1b\sigma_g(e_1)L_g=e_1b(p-e_1)L_g=be_1(p-e_1)L_g=0\]and
\[(p-e_1)bL_g(p-e_1)=(p-e_1)b\sigma_g(p-e_1)L_g=(p-e_1)be_1L_g=(p-e_1)e_1bL_g=0.
\]
If $p=I$, then $e_1$ and $I-e_1$ satisfy the conditions of lemma. If $p\neq I$, by lemma \ref{outer automorphism}, we have $\sigma_g|_{(I-p)\N (I-p)}\neq id$. By the above observations, there exists a unitary operator $u \in (I-p) \N (I-p)$ with $u=u^*$ such that \[\sigma_g(u)=-u.\] Let \[u=2e_2-(I-p).\] Since $\sigma_g(u)=-u$, we have \[\sigma_g(e_2)=(I-p)-e_2\] and \[\sigma_g((I-p)-e_2)=e_2.\] Note that $e_2 \leq I-p$ and $(I-p)b=b(I-p)=0$. Let \[e=e_1+e_2\] and \[f=p-e_1+(I-p)-e_2.\] Then $e+f=I$ and $\sigma_g(e)=f$. Thus
\begin{eqnarray}\nonumber
  ebL_ge &=& (e_1+e_2)bL_g(e_1+e_2) \\\nonumber
   &=& e_1bL_ge_1+e_1bL_ge_2+e_2bL_g(e_1+e_2) \\\nonumber
   &=& e_1b\sigma_g(e_1)L_g+e_1b\sigma_g(e_2)L_g \\\nonumber
   &=& e_1b(p-e_1)L_g+e_1b((I-p)-e_2)L_g\\\nonumber
   &=&0
\end{eqnarray}
and
\begin{eqnarray}\nonumber
  fbL_gf &=& ((p-e_1)+(I-p-e_2))bL_g((p-e_1)+(I-p-e_2)) \\\nonumber
   &=& (p-e_1)bL_g(p-e_1)+(p-e_1)bL_g(I-p-e_2) \\\nonumber
   &~& +(I-p-e_2)bL_g((p-e_1)+(I-p-e_2)) \\\nonumber
   &=& (p-e_1)b\sigma_g(p-e_1)L_g+(p-e_1)b\sigma_g(I-p-e_2)L_g\\\nonumber
   &=& (p-e_1)be_1L_g+(p-e_1)be_2L_g\\ \nonumber
   &=&0.
\end{eqnarray}
This completes the proof.

\end{proof}

\begin{Lemma}\label{adjoing operator can be cutted to zero}
Let $\M=\N\rtimes \mathbb{Z}_{2}$ and $a=a^*,b=b^*$ in $\N$ with $\sigma_g(a)=a,\sigma_g(b)=-b$. Then there is a family of finitely many mutually orthogonal nonzero projections $\{e_i \}_{i=1} ^{N}$ in $\mathcal{N}$ such that $\sum_{i=1}^{N}e_i =I$ and $e_i(a+ib)L_ge_i=0$ for $1\leq i\leq N$.
\end{Lemma}
\begin{proof}
Let $\lambda=\tau(a)$. Since $\tau(a-\lambda I)=0$ and the fixed point algebra $\N_1$ of $\sigma_{g}$ is a type $\rm II_1$ factor, by lemma \ref{WFY3}, there is a family of finitely many mutually orthogonal nonzero projections $\{E_j \}_{j=1}^{N}$ in $\N_1$ with $\sum_{j=1}^{N}E_j =I$ such that
\[E_j(a-\lambda I)E_j=0\] for $1\leq j\leq N$.
Then \[E_j(a-\lambda I)L_gE_j=E_j(a-\lambda I)\sigma_g(E_j)L_g=E_j(a-\lambda I)E_jL_g=0,~~ 1\leq j\leq N.\] Consider $E_j \N E_j$. Since $\sigma_g(E_jbE_j)=-E_jbE_j$, by lemma \ref{sigma b=-b}, we have two nonzero projections $e_{j}$ and $f_j$ in $E_j \N E_j$ with $e_j+f_j=E_j$ and $\sigma_g(e_j)=f_j$, such that \[e_jbL_{g}e_j=f_jbL_{g}f_j=0,~~j=1,\ldots,N.\]
At the same time, we have that $\sum_{j=1}^N (e_j+f_j)=\sum_{j=1}^NE_j=I$ and
\[e_j(ib+\lambda I)L_{g}e_j=f_j(ib+\lambda I) L_{g}f_j=\lambda f_j\sigma_g(f_j)L_g=\lambda f_je_jL_g=0,~~j=1,\ldots,N.\] Thus, for $1\leq j\leq N$, we have
\begin{eqnarray}\nonumber
  e_j(a+ib)L_ge_j &=& e_j((a-\lambda I)+(ib+\lambda I))L_g e_j \\\nonumber
   &=& e_j(a-\lambda I)L_g e_j+e_j(ib+\lambda I))L_g e_j \\\nonumber
   &=& e_j(a-\lambda I)f_jL_g \\\nonumber
   &=&e_j E_j(a-\lambda I)E_j f_jL_g\\ \nonumber
   &=&0
\end{eqnarray}
and
\begin{eqnarray}\nonumber
 f_j(a+ib)L_gf_j &=& f_j((a-\lambda I)+(ib+\lambda I))L_g f_j \\\nonumber
   &=& f_j(a-\lambda I)L_g f_j+f_j(ib+\lambda I))L_g f_j \\\nonumber
   &=& f_j(a-\lambda I)e_jL_g \\\nonumber
   &=&f_j E_j(a-\lambda I)E_j e_jL_g\\ \nonumber
   &=&0.
\end{eqnarray}
\end{proof}

\begin{Lemma}\label{conditional expectation to N}
Let $\N\subseteq \M$ be an inclusion of type $\rm II_{1}$ factors and $[\M:\N]=2$. For any $x\in\M$, $x=x^{*}$ with $E_{\N}(x)=0$, there is a family of finitely many mutually orthogonal nonzero projections $\{e_i \}_{i=1}^{N}$ in $\mathcal{N}$ such that $\sum_{i=1}^{N}e_i =I$ and $e_ixe_i=0$ for $1\leq i\leq N$.
\end{Lemma}
\begin{proof}
Since $E_{\N}(x)=0$, we have $x=z L_{g}$, where $z\in \N$. Let $z=a+ib$, where $a$ and $b$ are self-adjoint operators in $\N$. Since \[x^{*}=L_{g}z^{*}=\sigma_{g}(z^{*})L_{g}=\sigma_{g}(a-ib)L_{g}\] and  $x=x^{*}$, then \[(a+ib)L_{g}=\sigma_{g}(a-ib)L_{g}=(\sigma_{g}(a)-i\sigma_{g}(b)) L_{g}.\] This implies that \[\sigma_{g}(a)=a ~~\text{and}~~ \sigma_{g}(b)=-b.\] By lemma \ref{adjoing operator can be cutted to zero}, there is a family of finitely many mutually orthogonal nonzero projections $\{e_i \}_{i=1} ^{N}$ in $\mathcal{N}$ such that $\sum_{i=1}^{N}e_i =I$ and $e_ixe_i=0$ for $1\leq i\leq N$.

\end{proof}

\begin{Lemma}\label{conditional expectation to N1}
Let $\M=\N\rtimes \mathbb{Z}_{2}$, $\N_{1}$ be the fixed point algebra of $\sigma_{g}$ and $E_{\N_{1}}:\N\rightarrow\N_{1}$ be the conditional expectation. Then $\sigma_{g}\circ E_{\N_{1}}=E_{\N_{1}}\circ \sigma_{g}=E_{\N_{1}}$.
\end{Lemma}
\begin{proof}
Given $x\in \N$, the conditional expectation $E_{\N_{1}}$ is determined by the equation
\begin{equation}\label{equation 1}
\langle x-E_{\N_{1}}(x), y\rangle=0
\end{equation}
for all $y\in \N_{1}$. Note that $\sigma_{g}$ is a $*$-automorphism preserving the inner product. Applying $\sigma_{g}$ to both sides of equation (\ref{equation 1}) gives \[\langle \sigma_{g}(x)-\sigma_{g}(E_{\N_{1}}(x)), \sigma_{g}(y)\rangle=0\] for all $y\in \N_{1}$. Since $\N_{1}$ the fixed point algebra of $\sigma_{g}$, then \[\langle \sigma_{g}(x)-E_{\N_{1}}(x), y\rangle=0\] for all $y\in \N_{1}$. Hence $\sigma_{g}\circ E_{\N_{1}}(x)=E_{\N_{1}}(x)=E_{\N_{1}}\circ \sigma_{g}(x), \forall x \in \N$.

\end{proof}

\begin{Lemma}\label{conditional expectation to N1 0}
Let $\M=\N\rtimes \mathbb{Z}_{2}$, $\N_{1}$ be the fixed point algebra of $\sigma_{g}$ and $E_{\N_{1}}:\N\rightarrow\N_{1}$ be the conditional expectation. If $x\in \N$ and $\sigma_{g}(x)=-x$, then $E_{\N_{1}}(x)=0$.
\end{Lemma}
\begin{proof}
By lemma \ref{conditional expectation to N1}, we have
\[E_{\N_{1}}(x)=E_{\N_{1}}(\sigma_{g}(x))=E_{\N_{1}}(-x).\]
Then $E_{\N_{1}}(x)=0$.

\end{proof}

\begin{Lemma}\label{ei(1+x)Lgei=0}
Let $\M=\N\rtimes \mathbb{Z}_{2}$,  $\N_{1}$ be the fixed point algebra of $\sigma_{g}$ and $E_{\N_{1}}:\N\rightarrow\N_{1}$ be the conditional expectation. If $x\in \N$ and $\sigma_{g}(x)=-x$, then there is a family of finitely many mutually orthogonal nonzero projections $\{e_i \}_{i=1} ^{N}$ in $\mathcal{N}$ such that $\sum_{i=1}^{N}e_i =I$ and $e_{i}(\lambda+x)L_{g}e_{i}=0$ for any $\lambda\in\mathbb{C}$ and $1\leq i\leq N$.
\end{Lemma}
\begin{proof}
Let $x=a+ib$, where $a^*=a$, $b^*=b\in \N$. By lemma \ref{conditional expectation to N1 0}, we have $E_{\N_{1}}(x)=0$. Then $E_{\N_{1}}(b)=0$. Note that \[[\N:\N_{1}]=2.\] By lemma \ref{conditional expectation to N}, there is a family of finitely many mutually orthogonal nonzero projections $\{f_j \}_{j=1}^{N}$ in $\N_{1}$ such that $\sum_{j=1}^{N}f_j =I$ and \[f_{j}bf_{j}=0\] for $1\leq j\leq N$. Since $\sigma_{g}(f_{j})=f_{j}$ for $1\leq j\leq N$, we have
\[f_{j}(\lambda+x)L_{g}f_{j}=f_{j}(\lambda+x)\sigma_{g}(f_{j})L_{g}=f_{j}(\lambda+x)f_{j}L_{g}=
f_{j}(\lambda+a)f_{j}L_{g}.\]Note that $\sigma_g(a)=-a$, hence
\[\sigma_{g}(f_{j}af_{j})=-f_{j}af_{j}.\]
Then, by lemma \ref{sigma b=-b}, there are projections $e_{j1}, e_{j2}\in f_j\N f_j$ with $e_{j1}+e_{j2}=f_j$, $\sigma_g(e_{j1})=e_{j2}$ such that for all $1 \leq j \leq N$, \[e_{j1}(f_jaf_j)L_{g}e_{j1}=0\]
and
\[e_{j2}(f_jaf_j)L_{g}e_{j2}=0.\]
At the same time, we have for all $1 \leq j \leq N$,
\[e_{j1}(f_j(\lambda+a)f_j)L_{g}e_{j1}=0\]
and
\[e_{j2}(f_j(\lambda+a)f_j)L_{g}e_{j2}=0.\]
Let $e_k=e_{ji},1\leq j\leq N, i=1,2$. Therefore, we have $\sum_{k=1}^{2N}e_k=I$ and $e_k(\lambda+x)L_ge_k=0$ for any $\lambda \in \mathbb{C}$ and $1 \leq k \leq 2N.$
\end{proof}

\begin{Lemma}\label{index 2 finite proj in N cut to 0 single}
Let $\N\subseteq \M$ be an inclusion of type $\rm II_{1}$ factors and $[\M:\N]=2$. For any $x\in\M$ with $E_{\N}(x)=0$, there is a family of finitely many mutually orthogonal nonzero projections $\{e_i \}_{i=1} ^{N}$ in $\mathcal{N}$ such that $\sum_{i=1}^{N}e_i =I$ and $e_ixe_i=0$ for $1\leq i\leq N$.
\end{Lemma}
\begin{proof}
Since $E_{\N}(x)=0$, $x=zL_g$, where $z\in \N$. We decompose $x$ as:
\[z=\frac{z+\sigma_{g}(z)}{2}+\frac{z-\sigma_{g}(z)}{2}.\]
Denote $a=\frac{z+\sigma_{g}(z)}{2},b=\frac{z-\sigma_{g}(z)}{2}$. Then \[\sigma_g(a)=a ~~\text{and}~~ \sigma_g(b)=-b.\] Let $\lambda=\tau(a)$. By lemma \ref{WFY3}, there exists a family of finitely many mutually orthogonal nonzero projections $\{f_j \}_{j=1}^{N}$ in $\N_1$, where $\N_1$ is the fixed point algebra of $\sigma_g$, such that $\sum_{j=1}^{N}f_j =I$ and \[f_{j}(a-\lambda)f_{j}=0\] for $1\leq j\leq N$. Note that $\sigma_g(f_j)=f_j$, hence \[\sigma_g({e_jbe_j})=-e_jbe_j.\]
Then by lemma \ref{ei(1+x)Lgei=0}, there is a family of finitely many mutually orthogonal nonzero projections $\{e_{jk} \}_{k=1}^{N_j}$ in $\N$ such that $\sum_{k=1}^{N_j}e_{jk} =f_{j}$ and $e_{jk} (\lambda +b)L_ge_{jk}=0$ for $1\leq j\leq N$ and $1\leq k\leq N_j$. Therefore we have $\sum_{jk}e_{jk}=I$ and
\begin{eqnarray}\nonumber
e_{jk}x e_{jk}&=&e_{jk}((a-\lambda)+(\lambda+b))L_g e_{jk}=e_{jk}(a-\lambda)L_g e_{jk}\\ \nonumber
&=&e_{jk}f_j(a-\lambda)L_g f_j e_{jk}=e_{jk}f_j(a-\lambda)f_jL_g e_{jk}\\ \nonumber
&=& 0. \nonumber
\end{eqnarray}
Relabel $\{e_{jk}\}$ as $\{e_i\}_{i=1}^N$. Then $e_i\in \N$, $\sum_{i=1}^Ne_i=I$ and $e_ixe_i=0$ for $1\leq i\leq N$.
\end{proof}

\begin{Lemma}\label{index 2 finite proj in N cut to 0 many}
Let $\N\subseteq \M$ be an inclusion of type $\rm II_{1}$ factors and $[\M:\N]=2$. For any $x_{1},\ldots, x_{n}\in \M$ with $E_{\N}(x_j)=0,~j=1,\ldots,n$, there is a family of finitely many mutually orthogonal nonzero projections $\{e_i \}_{i=1} ^{N}$ in $\mathcal{N}$ such that $\sum_{i=1}^{N}e_i =I$ and $e_ix_je_i=0$ for $1\leq i\leq N$ and $1\leq j\leq n$.
\end{Lemma}
\begin{proof}
For $n=1$, by lemma \ref{index 2 finite proj in N cut to 0 single}, this lemma is true. Suppose that the lemma holds for any $k$ elements with $E_{\N}(x_j)=0$ for $j=1,\ldots,k$.

Now consider $k+1$ elements with $E_{\N}(x_j)=0$ for $j=1,\ldots,k+1$. By the inductive hypothesis, there exists a family of finitely many mutually orthogonal nonzero projections $\{f_i\}_{i=1}^{M}$ in $\mathcal{N}$ such that $\sum_{i=1}^{M}f_i =I$ and $f_ix_{j}f_i=0$ for $1\leq i\leq M$ and $1\leq j\leq k$.

For each $f_{i}$, we consider the element $f_{i}x_{k+1}f_{i}$. Since $E_{\N}(x_{k+1})=0$ and $f_i \in \N$, we have
\[E_{f_{i}\N f_{i}}(f_{i}x_{k+1}f_{i})=f_{i}E_{\N}(x_{k+1})f_{i}=0.\] Note that \[[f_{i}\M f_{i}:f_{i}\N f_{i}]=2.\]
By lemma \ref{index 2 finite proj in N cut to 0 single}, there exists a family of finitely many mutually orthogonal nonzero projections $\{e_{ij}\}_{j=1}^{L_{i}}$ in $f_{i}\N f_{i}$ such that $\sum_{j=1}^{L_{i}}e_{ij}=f_{i}$ and \[e_{ij}(f_{i}x_{k+1}f_{i})e_{ij}=0\] for $1\leq j\leq L_{i}$.
Note that
\[\sum_{i=1}^{M}\sum_{j=1}^{L_{i}}e_{ij}=\sum_{i=1}^{M}f_{i}=I.\]
Moreover, for each $e_{ij}$, we have
\[e_{ij}x_{l}e_{ij}=e_{ij}f_{i}x_{l}f_{i}e_{ij}=0
\]
for $1\leq l\leq k$, and
\[e_{ij}x_{k+1}e_{ij}=e_{ij}f_{i}x_{k+1}f_{i}e_{ij}=0.\]
Thus the family $\{e_{ij}\}$ satisfies the required conditions for $k+1$ elements. By induction, we have this lemma.
\end{proof}

\noindent\textbf{The proof of Theorem \ref{index 2 finite proj in N cut to 0 many trace}:}
Set \[x_{j}=E_{\N}(x_{j})+(x_j-E_{\N}(x_{j}))=y_{j}+z_{j},~j=1,\ldots,m.\]
Then $y_{j}\in \N$ and \[E_{\N}(z_{j})=0,~j=1,\ldots,m.\] Note that \[\tau(y_{j})=\tau(E_{\N}(x_{j}))=\tau(x_{j})=0,~j=1,\ldots,m.\]
By lemma \ref{WFY3}, there is a family of finitely many mutually orthogonal nonzero projections $\{p_i \}_{i=1} ^{N}$ in $\N$ such that $\sum_{i=1} ^{N}p_i =I$ and \[p_{i}y_{j}p_{i}=0\] for $1\leq i\leq N$ and $1\leq j \leq m$.
At the same time, \[E_{p_{i}\N p_{i}}(p_{i}z_{j}p_{i})=p_iE_\N (z_{j})p_{i}=0\]
and \[[p_{i}\M p_{i}:p_{i}\N p_{i}]=2.\] By lemma \ref{index 2 finite proj in N cut to 0 many}, there is a family of finitely many mutually orthogonal nonzero projections $\{e_{ik} \}_{k=1} ^{N_i}$ in $\mathcal{N}$ such that $\sum_{k=1}^{N_i}e_{ik} =p_{i}$ and \[e_{ik}z_{j}e_{ik}=0\] for $1\leq j\leq m$, $1 \leq i \leq N$ and $1 \leq k \leq N_i$. Thus the family $\{e_{ik}\}$ satisfies the required conditions of this theorem.
\qed

\vskip 0.5cm

In the following, we provide an application of Theorem 1.3. Let $\A$ be a Banach algebra with the identity. Given $a, b\in \A$, the Rosenblum operator $\tau_{ab}\in L(\A)$ is defined by \[\tau_{ab}(c)=ac-cb.\]

\begin{Lemma}\label{tau(ab)}{\rm (}\cite{D. A. Herrero} Cor 3.2{\rm )}
$\sigma(\tau_{ab})\subset \sigma(a)-\sigma(b)$.
\end{Lemma}

The following lemma is motivated by Lemma 2.2 of \cite{Marcoux}. We provide a complete proof for the reader¡¯s convenience.

\begin{Lemma}\label{diag operator is a commutator then oper is comm}
Given a pair of unital $C^{*}$-algebras $\N\subseteq \M$ and suppose that $p_{1}, p_{2},\ldots , p_{n}$ are mutually orthogonal projections in $\N$ with $p_{1}+p_{2}+\cdots+p_{n}=1$. For any $b\in \M$, if there exist $x_{k}\in p_{k}\N p_{k}$ and $y_{k}\in p_{k}\M p_{k}$ such that $p_{k}bp_{k}=[x_{k}, y_{k}]$, then there are $d\in \N$ and $z\in \M$ such that $b=[d, z]$.
\end{Lemma}
\begin{proof}
If $x_{k}\neq0$, then we can replace $x_{k}$ by $\frac{x_{k}}{\|x_{k}\|}$ and $y_{k}$ by $\|x_{k}\|y_{k}$, which allows us to assume that $\|x_{k}\| = 1$. In particular, we have that \[\sigma_{p_{k}Np_{k}}(x_{k})\subseteq\bar{\mathbb{D}}.\]
Let $d_{k}=x_{k}+3kp_{k}$, so that $d_{k}\in p_{k}\N p_{k}$ and \[\sigma_{p_{k}\M p_{k}}(x_{k})\subseteq\sigma_{p_{k}\N p_{k}}(d_{k})\subseteq 3k+\bar{\mathbb{D}}.\]
At the same time, we still have that \[[d_{k}, y_{k}]=[x_{k}+3kp_{k}, y_{k}]=[x_{k}, y_{k}]=p_{k}bp_{k}\] for all $k$.
If $j\neq k$, then \[\sigma_{p_{j}\M p_{j}}(d_{j})\cap\sigma_{p_{k}\M
p_{k}}(d_{k})=\emptyset.\] From this and lemma \ref{tau(ab)} it follows that for $1\leq j\neq k\leq n$, the Rosenblum operator
\[\begin{matrix}
\tau_{d_{j}, d_{k}}: & p_{j}\M p_{k} & \rightarrow  & p_{j}\M p_{k}\\
 & z & \mapsto  &d_{j}z-zd_{k}:=p_{j}bp_{k}
\end{matrix}\]
is invertible. Therefore, for each $1\leq j\neq k\leq n$, there exists $z_{jk} \in p_{j}\M p_{k}$ such that $\tau_{d_{j}, d_{k}}(z_{jk})= p_{j}bp_{k}$.
For $1\leq j \leq n$, let $z_{jj}= y_{j}$. Setting $d=\sum_{k=1}^{n}d_{k}$ and \[z=\sum_{1\leq j, k\leq n}z_{jk}=\sum_{1\leq j, k\leq n}p_{j}z_{jk}p_{k},\] then $d\in \N$, $z\in \M$, and it easily shows that
\begin{eqnarray}\nonumber
dz-zd
&=& \left(\sum_{j=1}^{n}d_{j}\right)\left(\sum_{1\leq j, k\leq n}z_{jk}\right)-\left(\sum_{1\leq j, k\leq n}z_{jk}\right)\left(\sum_{k=1}^{n}d_{k}\right) \\ \nonumber
&=& \sum_{1\leq j, k\leq n}d_{j}z_{jk}-z_{jk}d_{k}\\\nonumber
&=&  \sum_{1\leq j, k\leq n}p_{j}bp_{k}\\ \nonumber
&=& b.
\end{eqnarray}
\end{proof}

\begin{Corollary}
Let $\M=\N\rtimes \mathbb{Z}_{2}$ and $a\in \M$. If $\tau(a)=0$, then there are $d \in \N$ and $z \in \M$ such that $a=dz-zd$.
\end{Corollary}
\begin{proof}
Let $b=E_{\N}(a)\in \N$ and $a=b+(a-b):=b+c$. Then $\tau(b)=0$ and thus, by lemma \ref {WFY3}, there is a family of finitely many mutually orthogonal nonzero projections $\{p_i \}_{i=1} ^{n}$ in $\N$ such that $\sum_{i=1} ^{n}p_i =I$ and $p_{i}bp_{i}=0$ for $1\leq i\leq n$. At the same time,
\[
a=
\begin{pmatrix}
p_{1}cp_{1}&\ast  & \cdots&\ast \\
\ast& p_{2}cp_{2} & \cdots&\ast \\
\vdots&\vdots  &\ddots&\vdots \\
\ast& \ast & \cdots& p_{n}cp_{n}
\end{pmatrix}.
\]
Note that $E_{p_k\N p_k}(p_kcp_k)=0,k=1,\ldots,n$. Combining with lemma \ref{index 2 finite proj in N cut to 0 many} and lemma \ref{diag operator is a commutator then oper is comm}, we get this corollary.
\end{proof}

Combining with Theorem \ref{index 2 finite proj in N cut to 0 many trace} and Lemma \ref{diag operator is a commutator then oper is comm}, we have the following result.

\begin{Corollary}
Let $\N\subseteq \M$ be an inclusion of type $\rm II_{1}$ factors and $[\M:\N]=2$. For any $x_{1},\ldots, x_{n}\in \M$, if $\tau(x_j)=0,~j=1,\ldots,n$, then there are $d \in \N$ and $z_j \in \M$ such that $x_j=dz_j-z_jd$.
\end{Corollary}


\end{CJK}

\begin{thebibliography}{99}

\bibitem{AP} {C. Anantharaman and S. Popa}, An introduction to ${\rm II}_1$ factors, preprint.

\bibitem{RLA} R. Archbold, L. Robert and A. Tikuisis. The Dixmier property and tracial states for $C^{*}$-algebras. J. Funct. Anal., 273 (2017), no. 8, 2655-2718.

\bibitem{Alain Connes} A. Connes. Outer conjugacy classes of automorphisms of factors. Ann. Sci. \'{E}cole Norm. Sup., (4)8(1975), no. 3, 383-419.

\bibitem{Goldman} M. Goldman. On subfactors of factors of type $\rm II_1$. Mich. Math. J., 7(1960), 167-172.

\bibitem{HVG} H. Halpern, V. Kaftal and G. Weiss. The relative Dixmier property in discrete crossed products. J. Funct. Anal., 69 (1986), no. 1, 121-140.

\bibitem{D. A. Herrero} D. A. Herrero. Approximation of Hilbert space operators. I, Second ed., Pitman Res. Notes Math., vol. 224, Longman Sci. Tech., Harlow-New York 1989.


\bibitem{Jaeseong} J. Heo. Relative Dixmier property and complete boundedness of modular maps. Math. Z., 247(2004), no. 1, 89-99.



\bibitem{jones} V. F. R. Jones. Index for subfactors. Invent. Math., 72 (1983), no. 1, 1-25.

\bibitem{kR59} R. V. Kadison and I. M. Singer. Extensions of pure states. Amer. J. Math., 81(1959), 383-400.

\bibitem{kR} R. V. Kadison and J. Ringrose. Fundamentals of the theory of operator algebras, II. Pure Appl. Math. (N. Y.), Academic
Press, London, 1985.

\bibitem{Marcoux} L. W. Marcoux. Sums of small number of commutators. J. Operator Theory, 56(2006), no. 1, 111-142.

\bibitem{pop} F. Pop. Singular extensions of the trace and the relative Dixmier property in the type ${\rm II}_1$ factors. Proc. Amer. Math. Soc., 126(1998), no. 10, 2987-2992.

\bibitem{popa 81} S. Popa. On a problem of R. V. Kadison on maximal abelian *-subalgebras in factors. Invent. Math., 65(1981), 269-281.

\bibitem{popa 83} S. Popa. Singular maximal abelian subalgebras in continuous von Neumann algebras. J. Funct. Anal., 50(1983), 151-166.

\bibitem{popa 86} S. Popa. Correspondences. Preprint, 1986.

\bibitem{popa 90} S. Popa. Classification of subfactors: reduction to commuting squares. Invent. Math., 101(1990), 19-43.

\bibitem{popa 88} S. Popa. On amenable factors of type $\rm II_{1}$, in Operator Algebras and Applications. London Math. Soc. Lecture Note Ser., 136, pp. 73-184. Cambridge Univ. Press, 1988.

\bibitem{popa 94} S. Popa. Classification of amenable subfactors of type $\rm II$. Acta Math., 172(1994), no. 2, 163-255.

\bibitem{popa 95} S. Popa. Classification of subfactors and their endomorphisms, CBMS Lecture Notes Series in Math. 86(1994), A. M. S., Providence, RI, 1995.

\bibitem{Popa 2000} S. Popa. On the relative Dixmier property for inclusions of $C^{*}$-algebras. J. Funct. Anal., 171(2000), no. 1, 139-154.

\bibitem{SS} A. Sinclair and R. Smith. Finite von Neumann algebras and masas. London Mathematical Society Lecture Note Series, vol. 351, Cambridge University Press, Cambridge, 2008.

\bibitem{SZ} S. Str$\check{a}$til$\check{a}$ and L. Zsid\'{o}. Lectures on von Neumann algebras. Abacus Press, Bucharest, 1979.

\bibitem{WC} S. Wen and L. Cao. On isotropic subspaces of operators in ${\rm II}_1$ factors. J. Math. Anal. Appl., 538(2024), no. 1, Paper No.128332, 11 pp.

\bibitem{WFY3} S. Wen, J. Fang and Z. Yao. A stronger version of Dixmier's averaging theorem and some applications. J. Funct. Anal., 287(2024), no. 8, Paper No.110569, 13 pp.


\end{thebibliography}
\end{document}